\documentclass[english,american]{article}

\usepackage[letterpaper]{geometry}
\usepackage{float}
\usepackage{units}
\usepackage{amsthm}
\usepackage{amsmath}
\usepackage{graphicx}
\usepackage{amssymb}
\usepackage{datetime}
\usepackage[colorinlistoftodos, textwidth=1.5\marginparwidth, shadow]{todonotes}

\makeatletter

\theoremstyle{plain}
\newtheorem{thm}{Theorem}
\theoremstyle{definition}
\newtheorem{defn}[thm]{Definition}
\theoremstyle{plain}
\newtheorem{lem}[thm]{Lemma}
\theoremstyle{plain}
\newtheorem{prop}[thm]{Proposition}
 \ifx\proof\undefined\
   \newenvironment{proof}[1][\proofname]{\par
     \normalfont\topsep6\p@\@plus6\p@\relax
     \trivlist
     \itemindent\parindent
     \item[\hskip\labelsep
           \scshape
       #1]\ignorespaces
   }{%
     \endtrivlist\@endpefalse
   }
   \providecommand{\proofname}{Proof}
 \fi
\theoremstyle{plain}

\makeatother

\usepackage{babel}
\DeclareRobustCommand\transp{^{\mathrm{T}}}

\begin{document}

\title{\vspace*{-0.5in}
Random Triangle Theory with Geometry and Applications
}

\author{Alan Edelman and Gilbert Strang${\mbox{\hspace{-0.05in}}}$
\date{\today   \hspace{.1in} \currenttime}
\thanks{\noindent Department of Mathematics, Massachusetts Institute of Technology. \newline
{\hspace*{.2in} \it 2010 Mathematics Subject Classification.} 
Primary 51-xx, 15B52, 60B52, 52A22, 60G57, 51M15, 51M04. \newline 
{\hspace*{.2in} Keywords: random triangle, triangle space, obtuse triangle, Shape Theory, Hopf Fibration, Hopf Map, random matrix, Lewis Carroll, hemisphere, acute angle.} \newline
}}
\maketitle
\begin{abstract}
What is the probability that a random triangle is acute? We explore this old question from a modern viewpoint, taking into account linear algebra, shape theory, numerical analysis, random matrix theory, the Hopf fibration, and much much more. 
One of the best distributions of random triangles takes all six vertex coordinates as independent standard Gaussians. Six can be reduced to four by translation
of the center to $(0,0)$ or reformulation as a 2x2 random matrix problem.

In this note, we develop shape theory in its historical context for a wide audience.  We hope to encourage others to look again (and differently) at triangles.

We provide a new constructive proof, using the geometry of parallelians, of a  central result of shape theory:
{\bf  triangle shapes naturally fall on a hemisphere. }
We give several proofs of the key random result: that triangles are uniformly distributed
when the normal distribution is transferred to the hemisphere. A new proof  connects to the
distribution of random condition numbers.
Generalizing to higher dimensions, we obtain the ``square root ellipticity statistic" of random matrix theory.

 Another proof connects the Hopf map to the SVD of 2 by 2 matrices.
A new theorem describes  three similar triangles hidden in the hemisphere.
Many triangle properties are reformulated as matrix theorems, providing insight to both. 
 This paper argues for a shift of viewpoint 
to the modern approaches of random matrix theory.  As one example,  we
propose that the smallest singular value is an effective test for uniformity.
New software is developed and applications are proposed.

\pagebreak{}

\tableofcontents{}

\pagebreak{}
\end{abstract}

\section{Introduction}
Triangles live on a hemisphere and are linked to 2 by 2 matrices.
The familiar triangle is seen in a different light. New understanding
and new applications come from its connections to the modern developments of random matrix theory.
You may never look at a triangle the same way again.

We began with an idle
question: \textit{Are most triangles acute or obtuse\,?}
While looking for an answer, a note was passed in lecture. (We do
not condone our behavior\,!) The note contained an integral over a region in $\mathbb{R}^{6}$. The evaluation of that integral gave us a number -- the fraction of obtuse triangles. This paper will present several other ways to reach that number, but our real purpose is to provide a more complete picture of {}``triangle space.''

Later we learned that  Lewis Carroll (as Charles Dodgson) asked the same question in $1884$. His answer for the probability of an obtuse triangle (by his rules) was
\[
\dfrac{3}{8-\dfrac{6}{\pi}\sqrt{3}}\approx0.64.
\]
Variations of interpretation lead to multiple answers (see \cite{eisenberg96,Portnoy94}  and their references). Portnoy reports that in the first issue of The Educational Times $(1886)$, Woolhouse reached $9/8-4/\pi^{2}\approx0.72$. In every case obtuse triangles are the winners -- if our mental image of a typical triangle is acute, we are wrong. Probably a triangle taken randomly from a high school geometry book would indeed be acute. Humans generally think of acute triangles, indeed nearly equilateral triangles or right triangles, in our mental representations of a generic triangle. Carroll fell short of our favorite answer $3/4$, which is more mysterious than it seems. There is no paradox, just different choices of probability measure.

The most developed piece of the subject is humbly known as {}``Shape
Theory.'' It was the last interest of the first professor of mathematical statistics at Cambridge University, David Kendall
\cite{kendall89,kendall10}. We rediscovered on our own what the shape theorists knew, that \textit{triangles are naturally mapped onto points of a hemisphere}. It was a thrill to discover both the result and the history of shape space. 

We will add a purely geometrical derivation of the picture of triangle
space, delve into the linear algebra point of view, and connect triangles to
computational mathematics issues including condition number analysis, Kahan's \cite{kahan1,kahan2}  accuracy of areas of needle shaped triangles, and  random matrix theory.

We hope to rejuvenate the study of shape theory\,!
\begin{figure}[H]
\begin{centering}
\includegraphics[scale=0.25]{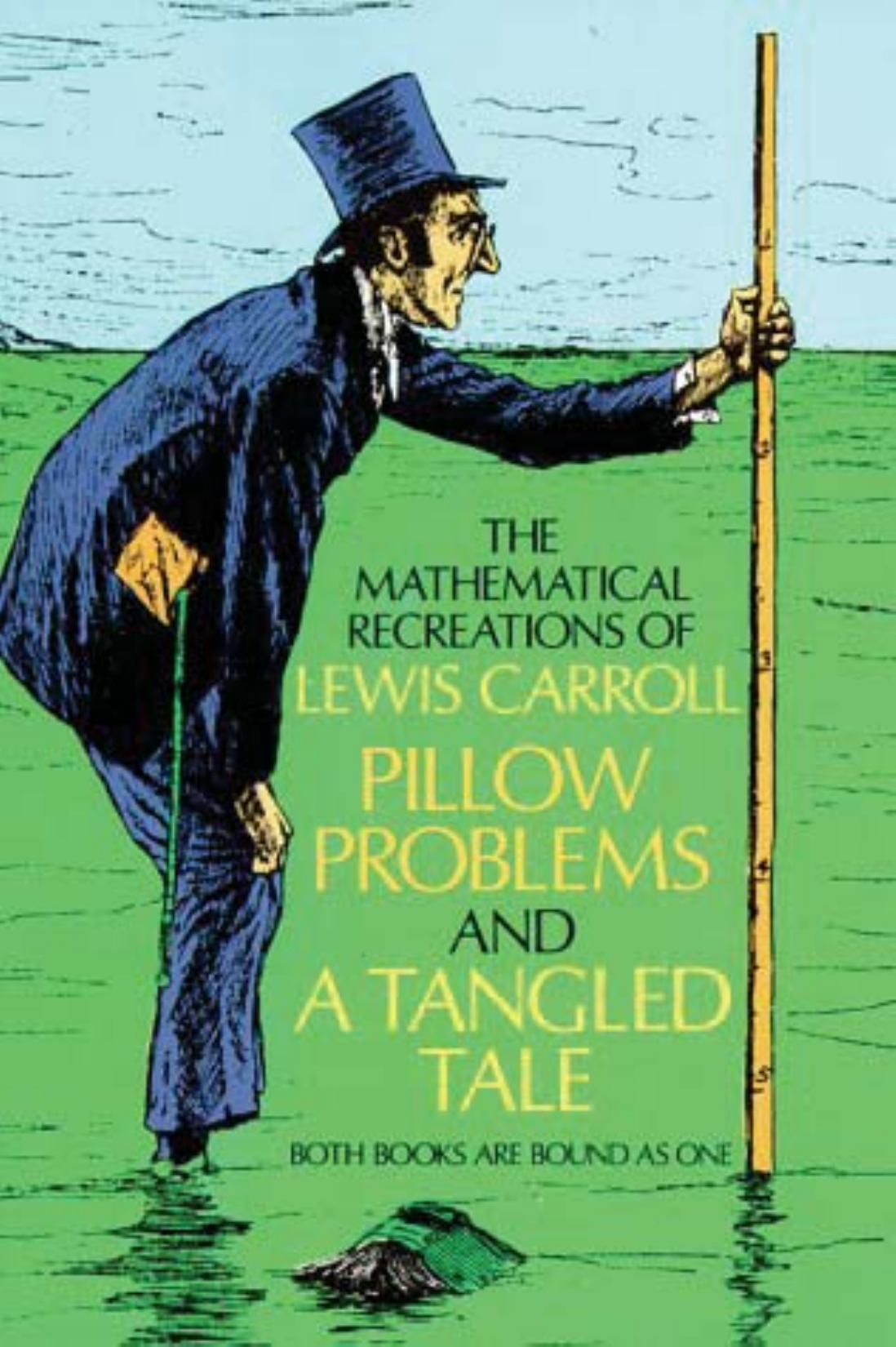}\includegraphics[angle=359,scale=0.22]{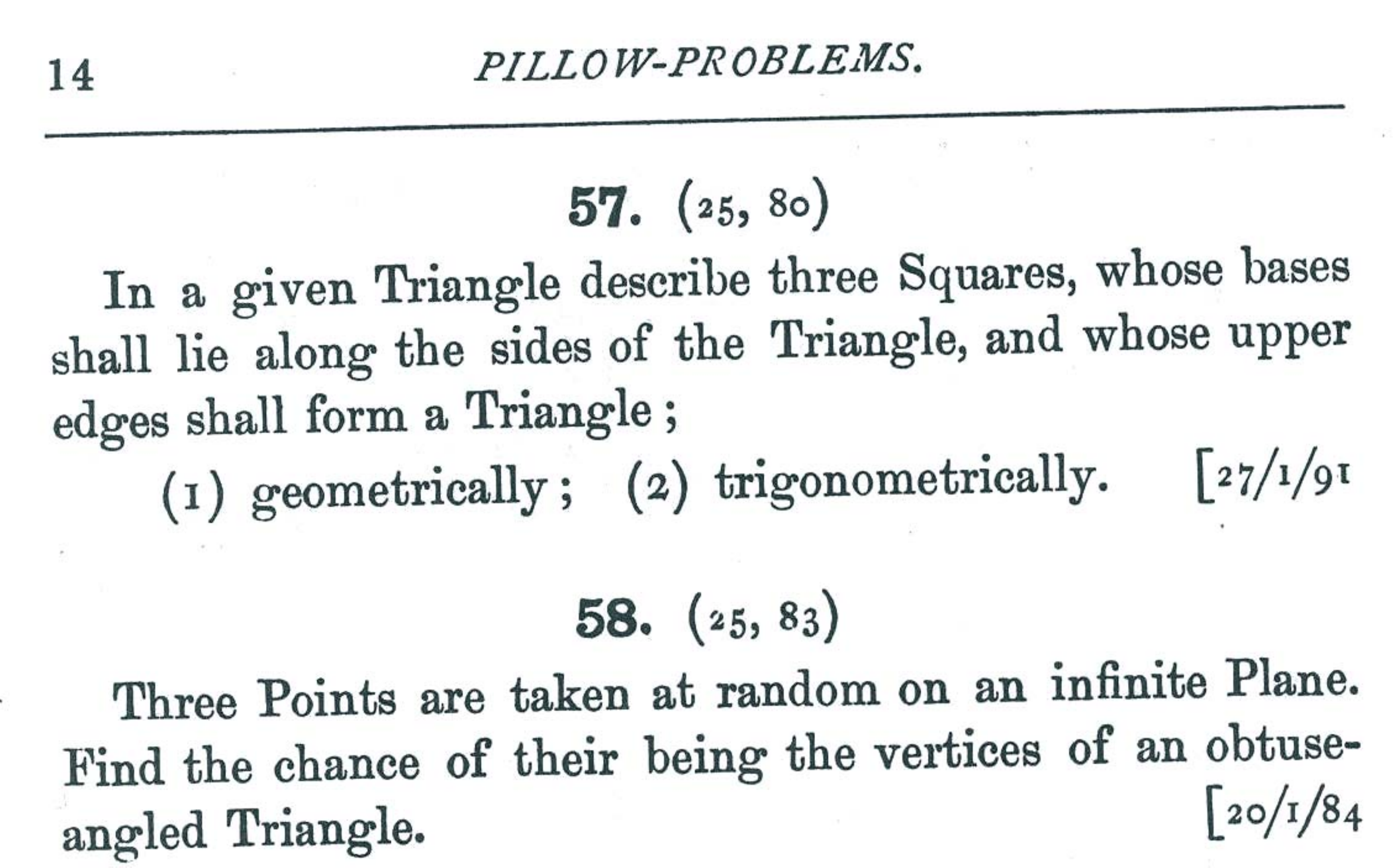}
\par\end{centering}
\caption[]{Lewis Carroll's Pillow Problem $58$ (January $20$, $1884$). $25$ and $83$ are page numbers for his answer and his method of solution. He specifies the longest side $AB$ and assumes that $C$ falls uniformly in the region where $AC$ and $BC$ are not longer than $AB$.}
\label{fig:Lewis-Carroll's-Pillow}
\end{figure}

\subsection{Random Angle Space}
\vspace{-.1in}

\begin{figure}[H]
\begin{center}
\includegraphics[scale=0.4]{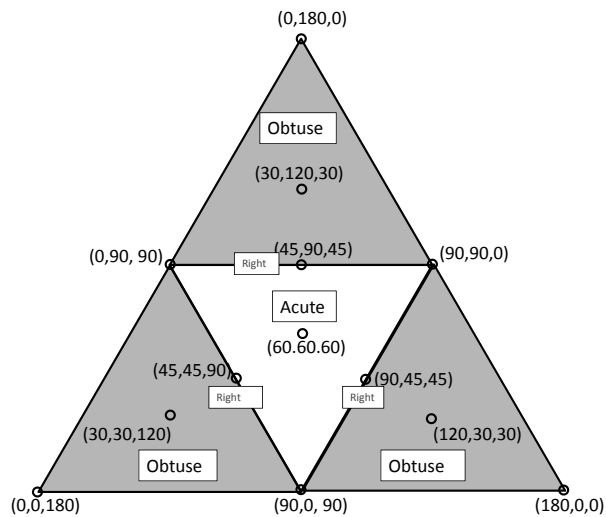}
\par\end{center}

\vspace{-2ex}
\caption[]{Uniform angle distribution\,: angle $1 +$
angle $2 +$ angle $3=180^{\circ}$.}
\label{fig:Uniform:-angle-1}
\end{figure}

The simplest  model  of a random triangle works with random angles. 
We mention it as a contrast to the Gaussian model that is more central to this paper.
Figure~\ref{fig:Uniform:-angle-1} shows a picture of angle space in $\mathbb{R}^{3}$. It is a {}``barycentric'' picture in the plane  for which $\hbox{angle }1+\hbox{ angle } 2+\hbox{angle }3=180^{\circ}$. A random point chosen from this space has a natural distribution, the uniform distribution, on the three angles. We then reach a simple fraction\,: $3/4$ of the triangles are obtuse.

\emph{The normal distribution on vertices also gives the fraction} $3/4$. This result is much less obvious. We are not aware of an argument that links the {}``angle picture'' with the normal distribution, though it is hard to imagine that anything in mathematics is a coincidence. Nonetheless, the normal distribution on vertices gives a very nonuniform distribution on angles.  We explore this further in Section 3.6.

In the next section we will discuss the {}``natural'' shape picture
of triangle space that is equivalent to the normal distribution. It seems worth repeating a key message: the uniform angle distribution is not the same triangle measure as the Gaussian distribution yet obtuse triangles  have the same probability, $3/4.$

\subsection{Table of Random Triangles}
\begin{figure}[H]
\begin{centering}
\includegraphics[scale=0.9]{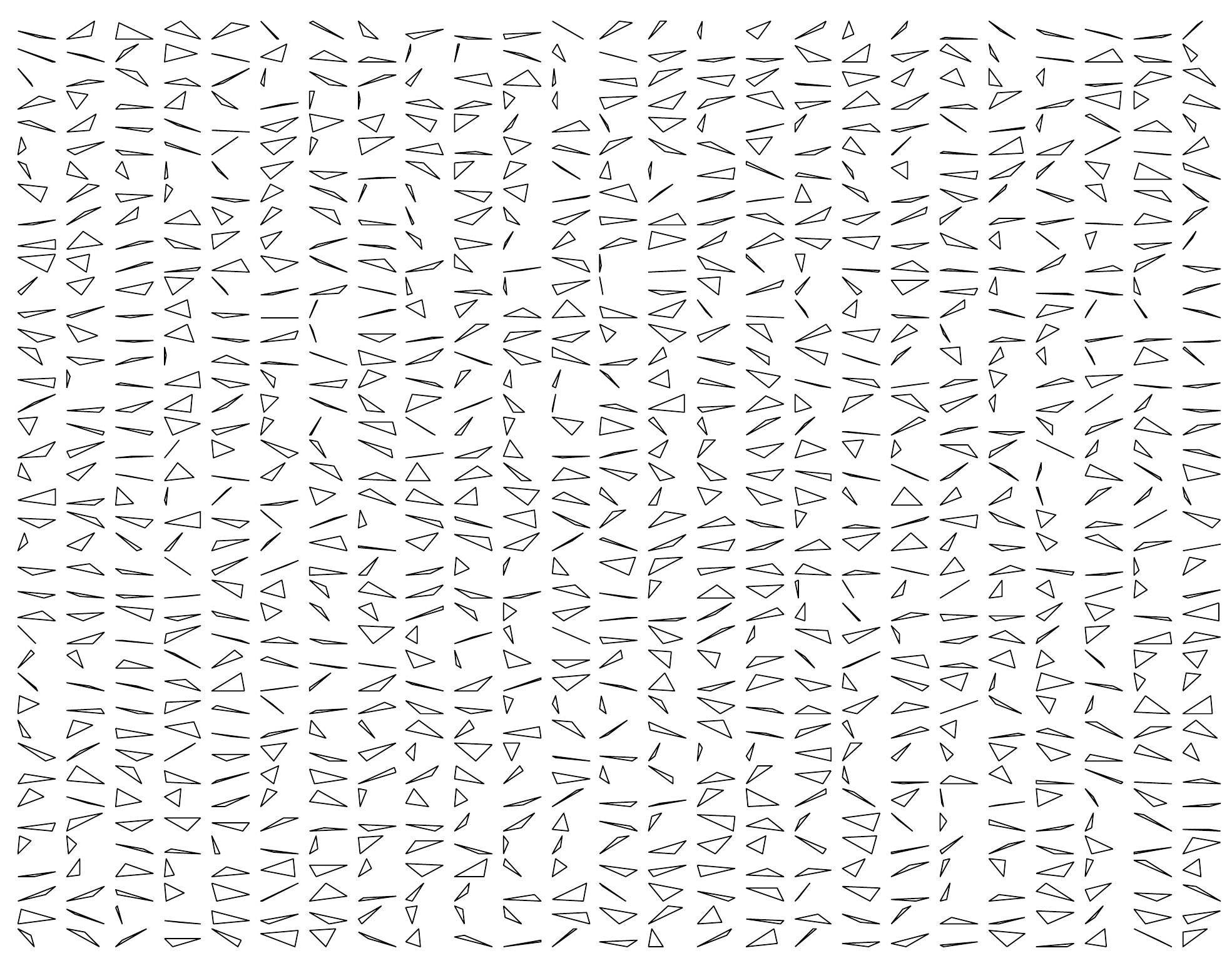}
\par\end{centering}
\caption[]{1000 Random Triangles (Gaussian Distribution): Most triangles are obtuse.}
\label{Fig:table}
\end{figure}

Before computers were handy, random number tables were widely available for applications. On first sight, a booklet of random numbers seemed an odd use of paper and ink, but these tables were highly useful. In the same spirit, we publish in Figure~\ref{Fig:table} a table of $1000$ random triangles. Only the shapes matter, not the scaling. You might try to count how many are acute. The vertices have six independent standard normally distributed coordinates. Each triangle is recentered and rescaled, but not rotated.

\subsection{A Fortunate {}``Optical Illusion''}
This subsection moves from $1,\!000$ to $50,\!000$ random triangles. The assumptions are the same\,: $x_{1},y_{1},x_{2},y_{2},x_{3},y_{3}$ are six independent random numbers drawn from the standard normal distribution (mean $0$, variance $1$). These six numbers generate a random triangle with vertices $(x_{1},y_{1})$, $(x_{2},y_{2})$ and $(x_{3,}y_{3})$. From the vertices, we can compute the three side lengths, $a$, $b$, $c$. And since we do not care about scaling, we may normalize  so that $a^{2}+b^{2}+c^{2}=1$.

\begin{figure}[H]
\begin{centering}
\includegraphics[scale=1.0]{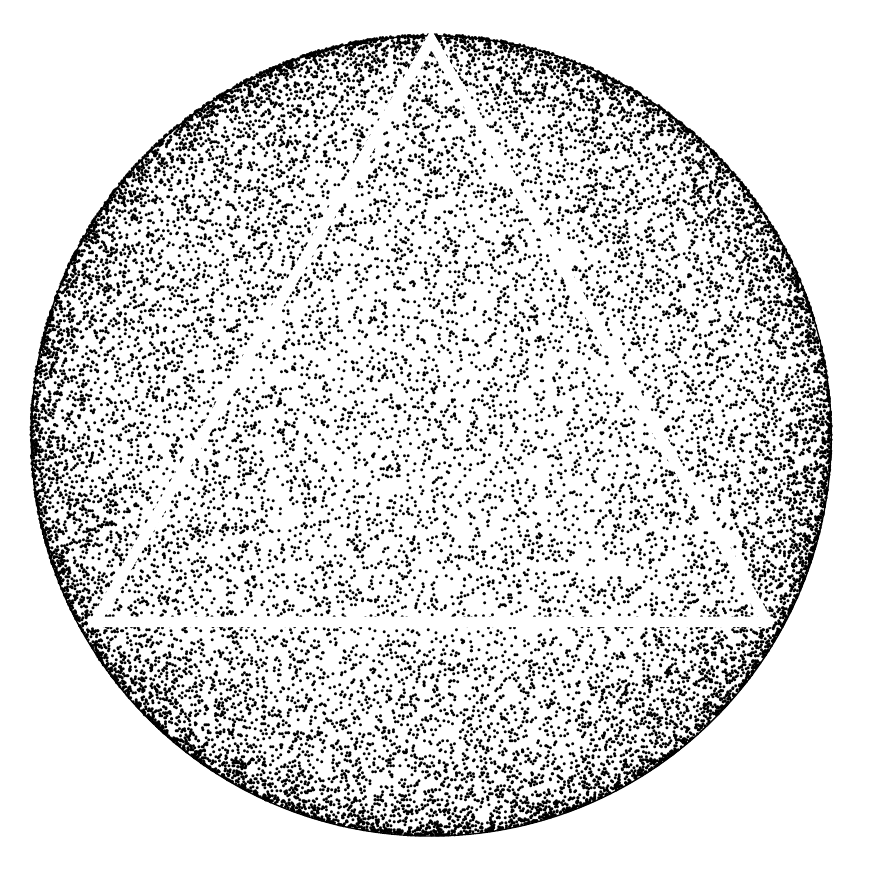}
\par\end{centering}
\caption[]{Points represent triangle shapes, where $(a^{2},b^{2},c^{2})$ is drawn in the plane $x+y+z=1$. Right triangles fall on  three white lines like $a^2+b^2=c^2=1/2$.} \label{fig:Points-represent-triangle}
\end{figure}

Instead of drawing the actual triangle, we represent it by
the point $(a^{2},b^{2},c^{2})$ in the plane $x+y+z=1$. The
result for many random triangles appears in Figure~\ref{fig:Points-represent-triangle}. 

The first curiosity is that the collection of points forms a disk.
This is the triangle inequality, though the connection  is not obvious. Points on the outer circumference
represent degenerate triangles, with area $0$. 

A second curiosity concerns right triangles. The points that represent
right triangles land on the white figure, an equilateral triangle
inscribed in the disk. 

A third curiosity, not visible in the picture, and also not obvious,  is that one quarter
of the points land in the acute region, inside the white equilateral
triangle. Each of the three disk segments representing obtuse triangles also contains one quarter of the points.

A fourth curiosity,  perhaps the most important of all, is 
the particular density of points (triangles) towards the perimeter of the disk. It is not difficult to imagine, in the spirit of many familiar optical illusions, that one is looking straight down towards the top of a hemisphere. The three white line segments are semicircles on the hemisphere, viewed {}``head on''.
Indeed, \textit{the points are uniformly distributed on a hemisphere}.
This fact was known to David Kendall, and is a major underpinning
of the subject of {}``shape theory.'' We discovered it ourselves
through this picture. Each area is $1/4$ of the area of the hemisphere, as Archimedes knew.

Grade school geometry emphasizes {}``side-side-side'' as enough
to represent any triangle, so what does height on the hemisphere represent\,? The answer is simple\,: height represents the area. The Equator has height zero, for degenerate triangles. The North Pole, an equilateral triangle, has maximum height. Latitudes represent triangles of equal area.

There is more. Triangles with one angle specified form small circles
on the hemisphere going through two vertices of the white figure. Indeed
it is good advice to take any triangle property and consider what
this hemisphere representation has to say.

\subsection{Two by Two Matrices: Turning Geometry into Linear Algebra}

A  triangle may be represented as a $2$ by $3$ matrix $T$ whose columns give the $x,y$ coordinates of the vertices. 
One way to remove two degrees of freedom is to translate the first vertex, say, to the origin.  A more symmetric way
to remove two degrees of freedom is to translate the centroid to the origin. 
There is something to be said for treating all vertices and also all
edges symmetrically. The classical law of sines does this, as we will do, but not
the law of cosines. Hero's formula for area  is also symmetric but   $\frac{1}{2}\hbox{ base }\times\hbox{ height }$ requires the choice of a base and hence is not symmetric.

Let $\Delta$ be the $2 \times 3$ matrix of  a reference equilateral triangle centered at the origin. 
This means that each of the three columns of $\Delta$ has the same euclidean length.
A $2\times2$ matrix $M$
transforms the triangle by taking the vertices of the equilateral triangle to the vertices of another triangle centered at the origin.
As illustrated  in
Figure~\ref{fig:2by2} below, every 
$2\times2$ matrix $M$ 
 may be associated with a
triangle with zero centroid 
through the equation $T=M\Delta$.  If the matrix $M$ is random, then the associated triangle $T$ is random.

\begin{figure}[H]
\centering
\includegraphics[scale=0.5]{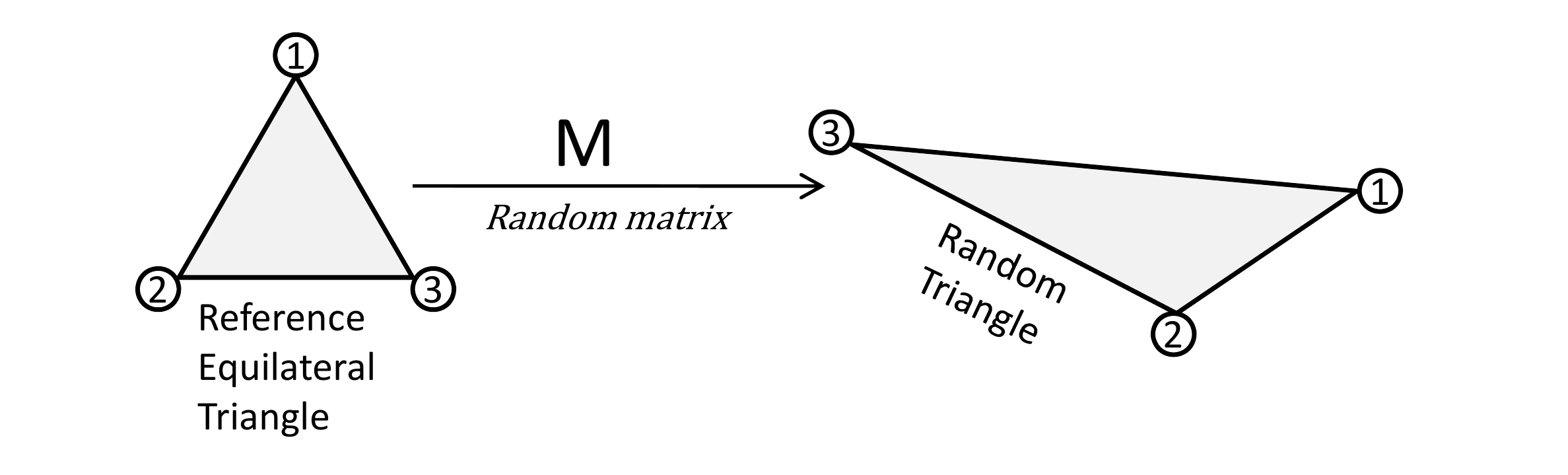}
\caption[]{Equivalence between $2\times 2$ matrices $M$ and triangles.}
\label{fig:2by2}
\end{figure}

An origin centered triangle  may alternatively be represented by its edge  vectors which
we may place in a $2 \times 3$ matrix $E$  whose columns add to the zero column  $(0,0)^T$.
We may again let $\Delta$ be the notation we use, this time for the edge vector matrix for an equilateral triangle.
Whether $\Delta$ is thought of as the vertices of an origin centered equilateral triangle or the edges, the important
property is that the columns of $\Delta$ are all of equal length. The vertices are equidistant from $(0,0)$, the edges have equal 
lengths.

If $M$ is a random $2 \times 2$ matrix, then $E=M\Delta$ produces a random triangle whose columns are triangle edges.
The columns  of $\Delta$ sum to $(0,0)^T$, hence the columns of $E$ sum to $(0,0)^T$, so $E$ encodes a proper closed triangle.
The edge lengths are the square roots of the diagonal elements of $E^{T}E$.  The vertex view would have produced $T$, whose three column vectors originate at the  centroid, giving less access to the edge lengths.

Section~2.2.1  takes a close look at  the choice of the Helmert matrix  $\Delta=
\left(\begin{array}{crc}
1/\sqrt{2} & -1/\sqrt{2} & 0\\
1/\sqrt{6} & 1/\sqrt{6} & -2/\sqrt{6}\end{array}\right)
$  
as our reference.

Consider as an example, the ``45,45,90"  zero centered  right  triangle with vertices $(-2,-1)$,$(1,-1)$,$(1,2)$
and edge lengths $3\sqrt{2},3,3$. 
 We form the  matrix $T=\left(\begin{array}{crc}
-2 & 1 & 1\\
-1& -1 & 2 \end{array}\right).
$
With the Helmert matrix  as $\Delta$, the corresponding $M_{v}=T\Delta^T$ for the vertex view is  $M_v=
-\left(\begin{array}{cc}
\sqrt{9/2} & \sqrt{3/2}   \\
0 & \sqrt{6}  \end{array}\right).
$
Any translation of the triangle produces the same $M_v$.

An edge view requires differences of columns of $T$.  
We may transform between 0 centered  vertices and edges with the equations
$$E = T 
\left(
\begin{array}{rrr}
1 & -1 & 0 \\
0 & 1 & -1 \\
-1 & 0 & 1
\end{array}
\right) \ \  \ {\rm and }
\ \  \ 
T = \frac{1}{3}E
\left(
\begin{array}{rrr}
1 &  0& -1\\
-1 & 1 & 0 \\
0 & -1  & 1
\end{array}
\right) 
. $$
These matrices are pseudoinverses as the $(1,1,1)$ direction is irrelevant. 

The $3\sqrt{2},3,3$ triangle is then  represented
as  $E=\left(\begin{array}{rrr}
 -3 & 3 & 0\\
  -3 & 0 & 3 \end{array}\right).
$
The corresponding edge view matrix is  $M_e=E\Delta^T=
- \left(\begin{array}{cc}
\sqrt{18}  &0  \\
\sqrt{9/2} & \sqrt{27/2}   
\end{array}\right).
$

It is possible to check that 
$M_v = M_e 
 \left(\begin{array}{rr}
1/2   &  \sqrt{3}/6    \\
-\sqrt{3}/6 & 1/2
\end{array}\right)
$
always holds.

This paper will take the edge view unless noted otherwise.

\section{Triangular Shapes $=$ Points on the Hemisphere}
Now we come to the heart of the paper. Every triangular shape with ordered vertices is naturally identified with a point on the hemisphere. {}``Random triangles are uniformly distributed on that hemisphere.'' We know of three constructions that exhibit the identification\,:
\begin{itemize}
\item Complex Numbers 
\item Linear Algebra\,: The Singular Value Decomposition
\item High School Geometry
\end{itemize}
\begin{flushleft}
The linear algebra approach, through the SVD, generalizes most readily
to higher dimensional shape theory. It also invites applications  using Random Matrix Theory which has advanced
considerably since the early inception of shape theory. (See Sections 3.2 and 4.1 of this paper for connections to Random Matrix Theory including random condition numbers and uniformity tests.)  Yet another benefit is the connection to the Hopf fibration.
\par\end{flushleft}

The geometric construction has two benefits not enjoyed by the others\,:

1. It can be understood with only knowledge of ordinary Euclidean
geometry 

2. It reveals the triangular shape and the point on the hemisphere in
the same picture.

\subsection{Complex Numbers}
A generic triangle can be scaled, rotated, and reflected so as to have vertices with complex coordinates $0,1,\zeta$. There are six possibilities for $\zeta$ in the upper half plane, corresponding to the six permutations of the vertices. One can consider $\bar{\zeta}$ in the lower half plane, to pick up six more possibilities.

The usual correspondence with the hemisphere is that the stereographic
projection maps $\zeta$ from the upper half plane to the hemisphere 
\cite{kendall89,kendall10}.

\subsection{Linear Algebra}
\subsubsection{The Helmert matrix and regular tetrahedra in $n$ dimensions}
The Helmert matrix $\Delta$ appears in statistics, group theory, and geodesy. It is a particular construction with orthonormal rows perpendicular to $(1,1,1)$\,:
\[
\Delta=
\left(\begin{array}{crc}
1/\sqrt{2} & -1/\sqrt{2} & 0\\
1/\sqrt{6} & 1/\sqrt{6} & -2/\sqrt{6}\end{array}\right).
\]
The columns contain the vertices of an equilateral triangle.  They are also the edges of a similar equilateral triangle. The operator view is that $\Delta^T$ is
an isometry (distances and angles are preserved) from
$\mathbb{R}^{2}$ to the plane $x+y+z=0$. This means that $\Delta\Delta^{T}=I_{2}$ and $\Delta^{T}\Delta=I_{3}-J_{3}/3$, where $J_{3}=$\texttt{ones(3)} is the all-ones matrix.

A useful property related to the duality between vertices and edges
is that
\[
\Delta\left(\begin{array}{rrr}
1 & -1\\
 & 1 & -1\\
-1 &  & 1\end{array}\right)=
\sqrt{3}\left(\begin{array}{cr}
\cos\frac{\pi}{6} & -\sin\frac{\pi}{6}\\
\sin\frac{\pi}{6} & \cos\frac{\pi}{6}\end{array}\right)\Delta.
\]
The left hand side takes the vector differences of
the columns of $\Delta$ (the 
vertices of
the equilateral triangle). These differences are the edges. We may place them as
vectors at a common origin. The right hand side is $\sqrt{3}$ times
a rotation matrix times $\Delta$.

The $\sqrt{3}$ represents  the familiar fact that
the edge lengths  of  an equilateral triangle are $\sqrt{3}$
times the distances of the vertices from the centroid.
Also if you rotate an edge counterclockwise by $\pi/6$, you  get the direction of a vertex.

The generalization of  Helmert's construction to $\mathbb{R}^{n}$
is the $n-1\times n$ matrix $\Delta_n$. It has  orthonormal rows, and the zeros
make it lower Hessenberg:
\begin{equation}
\Delta_{n}=\left(\begin{array}{cccc}
\sqrt{1\cdot2}\\
 & \sqrt{2\cdot3}\\
 &  & \ddots\\
 &  &  & \sqrt{(n-1)\cdot n}\end{array}\right)^{-1}\left(\begin{array}{crrcl}
1 & -1 & 0 & \ldots & 0\\
1 & 1 & -2 & \ldots & 0\\
\ldots & \ldots & \ldots & \ddots & \vdots\\
1 & 1 & 1 & \dots & -(n-1)\end{array}\right).
\label{eq:helmert}
\end{equation}
A transposed (and negated) form of that last matrix is available
in the statistics language R as \texttt{contr.helmert(n)}. The name indicates a matrix of {}``contrasts'' as used in statistics.

The square (and orthogonal) Helmert matrix adds a first row with entries $(1,\ldots,1)/\sqrt{n}$. It is obtained in MATLAB as an orthogonal test matrix with the command \texttt{gallery('orthog',n,4)}.

The rectangular $\Delta_{n}$ contains the vertices of a regular simplex from the {}``column view.'' When $n=4$ this is the regular tetrahedron with four vertices in $\mathbb{R}^{3}$. Our equilateral triangle $\Delta$ is $\Delta_{3}$. (The {}``row view'' consists of $n-1$ orthogonal rows each perpendicular to $(1,1,\ldots,1)$ in $\mathbb{R}^{n}$.) Again $\Delta_{n}\Delta_{n}^{T}=I_{n-1}$ and $\Delta_{n}^{T}\Delta_{n}=I_{n}-J_{n}/n$.

\subsubsection{The SVD\,: How $2\times2$ matrices connect triangles to the hemisphere}\label{sub:The-SVD:-How}
Let $M=U\Sigma V^{T}$ be the Singular Value Decomposition of a $2\times2$ matrix $M$. As shown in Figure~\ref{fig:2by2},  with
the edge viewpoint, the columns of $E=M\Delta$ are ordered edges of a triangle. We take a few steps to make the shape and the SVD unique\,:
\begin{itemize}
\item Scaling\,: We may assume that the squared entries of $M$ have sum $1$. The diagonal matrix $\Sigma$ then has  $1=\sigma_{1}^{2}+\sigma_{2}^{2}$.  The triangle edges then have sum of squares $1$ since ${\rm tr}(E^T\! E)={\rm tr}(\Delta^T M^T M\Delta)
=  {\rm tr}( M^T M\Delta\Delta^T) = {\rm tr}( M^T M)=1.$
\item 
The orthogonal factor $U$ in the SVD is unimportant, since $M$ and $U^{-1}M$ correspond to
the same triangle just
rotated or reflected. 
\item To make the SVD unique, assume that $\sigma_{1}\ge\sigma_{2}\ge0$ and that $V=\left(\begin{array}{rr}
\cos\theta & -\sin\theta\\ \sin\theta & \cos\theta\end{array}\right)$ with $0\le\theta<\pi$. There is a singularity in the SVD when $M=\frac{1}{\sqrt 2}I_{2}$. In this case $V$ is arbitrary.
\end{itemize}
We associate $V$ and $\Sigma$  with a point on the hemisphere of radius
$\frac{1}{2}$\,: The longitude is $2\theta$ and the height is $\sigma_{1}\sigma_{2}$. The latitude is thus asin($2\sigma_{1} \sigma_{2}$). From \foreignlanguage{english} {$1=\sigma_{1}^{2}+\sigma_{2}^{2}$} we have $0\le\sigma_{1}\sigma_{2}\le\frac{1}{2}$. The singularity of the North Pole (every angle is a longitude) is consistent with the singularity of the SVD at $\frac{1}{\sqrt 2}I_{2}$.

The height $\sigma_1\sigma_2$ is also $\left(\kappa+\kappa^{-1}\right)^{-1}$, where $\kappa=\sigma_{1}/\sigma_{2}\ge1$ is the condition number of $M$. The best conditioned matrix  $(\kappa=1)$ is $\frac{1}{\sqrt{2}}I_{2}$ at the North Pole, and corresponds to the equilateral triangle. The ill-conditioned matrices with $\kappa=\infty$, and at height $0$, are on the equator. They correspond to degenerate triangles, with  collinear vertices. We see our first hint of the  link between {}``ill-conditioned'' triangles and {}``ill-conditioned'' matrices.

One can now go either way from triangles to the hemisphere through the
$2\times2$ matrix $M$\,:
\begin{itemize}
\item Hemisphere to matrix to triangles\,: Start with a point on the hemisphere. Create $M=\Sigma V^{T}$ from the latitude asin($2\sigma_{1} \sigma_{2}$) and longitude $2\theta$. The edges of the triangle are the columns of $M\Delta$.
\item Triangles to matrix to hemisphere\,: Start with a triangle centered at $0$ whose squared edges sum to $1$. Then $M=$ ($2\times3$ matrix of triangle edges)$\Delta^{T}$. The point on the hemisphere comes from the SVD of $M$ by ignoring $U$, taking the height from $\det M=\sigma_{1}\sigma_{2},$ and the longitude from $\theta$.
\end{itemize}
The matrix $\mbox{\ensuremath{M}}$ is related to  the preshape in \cite{kendall89}.
 In that view, any triangle can be moved into a
standard position through $\Delta^{T}$. Our recommended view as shown in Figure~\ref{fig:2by2} is that the preshape is an operator mapping the equilateral triangle into a particular triangle. 
This viewpoint emphasizes the linear operator  that transforms triangles.

\subsection{Formulas at a glance}

As a summary, one may sequentially follow these steps to derive the
key formulas:

\begin{enumerate}
\item[] Reference Triangle:
edges from the three columns of \foreignlanguage{american}{ $\Delta=\left(\begin{array}{crc}
1/\sqrt{2} & -1/\sqrt{2} & 0\\
1/\sqrt{6} & 1/\sqrt{6} & -2/\sqrt{6}\end{array}\right).$}

\selectlanguage{american}%
Random Triangle: edges from the columns of $E=M\Delta.$ (Note $M=E \Delta^T$  since $\Delta \Delta^T=I$.) 

\item[{{(1})}] \selectlanguage{english} SVD  ($\sigma_1,\sigma_2,\theta$)  of Matrix $M$ ($\|M\|=1$) :

$M=\Sigma V^{T}=\left(\begin{array}{rr}
\sigma_{1}\\
 & \sigma_{2}\end{array}\right)\left(\begin{array}{rr}
\cos\theta & \sin\theta\\
-\sin\theta & \cos\theta\end{array}\right),$

$1 \ge  \sigma_{1}\ge\sigma_{2}\ge 0,\ \sigma_{1}^{2}+\sigma_{2}^{2}=1,\ 0\le\theta<\pi,$\foreignlanguage{american}{\vspace{0in}} ($U$ not needed.)

\item[{{(2})}] \selectlanguage{english} Triangle edges $a,b,c$:

$a^{2}+b^{2}+c^{2}=1,\ \left(\begin{array}{c}
a^{2}\\
b^{2}\\
c^{2}\end{array}\right)=\mbox{diag}(\Delta^{T}M^{T}M\Delta), $  $a+b\ge c$, $b+c \ge a$, $c+a \ge b$. \foreignlanguage{american}{\vspace{0in}}

\item[{{(3})}] \selectlanguage{english} Hemisphere of radius $1/2$ (coordinates $\lambda$, $\phi$ denote
$\frac{1}{2}(\cos \lambda \cos \phi,\cos \lambda \sin \phi, \sin \lambda $): 

Latitude = $\lambda=\mbox{asin}(2\sigma_{1}\sigma_{2}),$ Longitude $=\phi=2\theta,$
Height = $\frac{1}{2}\sin(\lambda)$. \foreignlanguage{american}{\vspace{0in}}
$0 \le \lambda \le \pi/2$,  $ \ 0 \le \phi <2 \pi$.

\item[{{(4})}] \selectlanguage{english} Disk of radius 1/2 (projection of hemisphere to polar coordinates $r,\phi$ ):

$r=\frac{1}{2}\cos(\lambda)$.  The angle $\phi=2\theta$ is the longitude of the point on the sphere. $0 \le r \le 1/2$, $0 \le  \phi  < 2 \pi$.
\end{enumerate}

The  area formulas for these descriptions:
\[
K=\mbox{Area = \ensuremath{\frac{\sigma_{1}\sigma_{2}}{\sqrt{12}}=\frac{1}{4}\sqrt{1-2(a^4+b^4+c^4)}=\sqrt{\frac{1-4r^{2}}{48}}=\frac{\mbox{Height}}{\sqrt{12}}=\frac{\mbox{sin($\lambda$)}}{\sqrt{48}}=(\kappa+\kappa^{-1})^{-1}/\sqrt{12}}.}\]

Here are conversion formulas. The (2)$\rightarrow$(4) and (4)$\rightarrow$(2) blocks are developed in Section 2.3.1.

\vspace{0.2in}
\hspace{-0.4in}\begin{tabular}{|l@{}|c@{}|c@{}|c@{}|c@{\,}|}
\hline 
From\textbackslash{}To & (1) SVD($M$) & (2) Triangle & (3) Hemisphere & (4) Disk\tabularnewline
\hline
\hline 
$\begin{array}{@{\,}l@{\,}}
\mbox{(1) SVD($M$)} \\
\hspace{0.2in} \sigma_{1},\sigma_{2},\theta\end{array}$ &  & $\begin{array}{@{\,}l@{\,}}
\rule{0pt}{0.15in} \ \ \ \ \ \ \ \ \mbox{diag}((M\Delta)^T(M\Delta)):\\[0.03in]
a^{2}=\frac{1}{3}(1-(\sigma_{1}^{2}-\sigma_{2}^{2})\cos2\theta_{+})\\[0.04in]
b^{2}=\frac{1}{3}(1-(\sigma_{1}^{2}-\sigma_{2}^{2})\cos2\theta_{-})\\[0.04in]
c^{2}=\frac{1}{3}(1-(\sigma_{1}^{2}-\sigma_{2}^{2})\cos2\theta)\\[0.05in] \end{array}$ & $\begin{array}{@{\,}l@{\,}}
\mbox{$\lambda=$asin(\ensuremath{2\sigma_{1}\sigma_{2})}}\\
\mbox{$\phi$=\ensuremath{2\theta}}\\[0.07in]
\mbox{Height=\ensuremath{K\sqrt{{12}}}}\\
\mbox{\ensuremath{=\sigma_{1}\sigma_{2}}=(\ensuremath{\kappa}+\ensuremath{\kappa^{-1})^{-1}}}\end{array}$ & $\begin{array}{@{\,}l@{\,}}
\rule{0pt}{0.15in} r\sin \phi=(M^T\! M)_{12}
\\
r\cos \phi=\frac{ (M^{\!T}\!\! M)_{11}-(M^{\!T}\!\! M)_{22} }{2}
\\ [0.05in]

r=\sqrt{\frac{1}{4}-\sigma_{1}^{2}\sigma_{2}^{2}}\\[.05in]
\ \ \ =\frac{1}{2}(\sigma_1^2-\sigma_2^2)\\[.03in]
\phi=2\theta\end{array}$\tabularnewline
\hline 
$\begin{array}{@{\,}l@{\,}}
\mbox{(2) Triangle}\\
\hspace{0.2in} a^{2},b^{2},c^{2}\end{array}$ & $\begin{array}{@{\hspace{-.1in}}l@{\,}}
\mbox{ 1)Triangle\ensuremath{\rightarrow}\ disk}\\
\mbox{ 2)then use (4)}\downarrow\end{array}$ & 
& $\! \! \frac{1}{\sqrt{6}}  \left( \begin{array}{@{\,}l@{\,}}
\cos \lambda \sin \phi \\ \cos \lambda \cos \phi \end{array}\right)=
\Delta\! \left( \begin{array}{@{\,}c@{\,}} {\rule{0pt}{0.15in} a^2} \\ b^2 \\c^2 \end{array}\right)$
&
$r\left( \begin{array}{@{\,}c@{\,}} \sin \phi \\ \cos \phi \end{array}\right)
=
\sqrt \frac{3}{2} \Delta \! \left( \begin{array}{@{\,}c@{\,}} {\rule{0pt}{0.15in} a^2} \\ b^2 \\c^2 \end{array}\right)
$ 
\vspace{0.02in}
\tabularnewline
\hline 
$\begin{array}{@{\,}l@{\,}}
\mbox{(3) Hemisphere}\\
\hspace{0.2in} \mbox{$\lambda,\phi$}\end{array}$ & 
$\begin{array}{@{\,}l@{\,}}
\sigma_1=\cos(\lambda/2) \\
\sigma_2=\sin(\lambda/2)
\end{array}$
& $\begin{array}{@{\,}l@{\,}}
\mbox{ 1) }r=\cos (\lambda)/2\\
\mbox{ 2) then use }\downarrow\end{array}$ 
& \rule{0pt}{0.2in} & $\begin{array}{@{\,}l@{\,}}
r=\cos(\lambda)/2\\
\phi=\phi\end{array}$\tabularnewline
\hline 
$\begin{array}{@{\,}l@{\,}}
\mbox{(4) Disk}\\
\hspace{0.2in} r,\phi\end{array}$ & $\begin{array}{@{\,}l@{\,}}
\rule{0pt}{0.14in} \sigma_{1}^{2}=\nicefrac{1}{2}+r\\
\sigma_{2}^{2}=\nicefrac{1}{2}-r\\
\theta=\phi/2\end{array}$ & $\begin{array}{@{\,}l@{\,}}
\rule{0pt}{0.14in} a^{2}=(1-2r\cos\phi_{+})/3\\
b^{2}=(1-2r\cos\phi_{-})/3\\
c^{2}=(1-2r\cos\phi)/3 
\\
\end{array}$ &
$\begin{array}{l}
\lambda=\mbox{acos}(2r) \\ \phi=\phi \end{array} $
 & \tabularnewline
\hline
\end{tabular}

$\theta_{\pm}=\theta\pm\pi/3,$ $\kappa=\sigma_{1}/\sigma_{2},$ $\phi_{\pm}=\phi\pm2\pi/3.$

\vspace{0.2in}

Also useful is the direct mapping of $M$ to the hemisphere in Cartesian coordinates:
$$
M=
\left(
\begin{array}{cc}
M_{11} & M_{12} \\
M_{21} & M_{22}
\end{array}
\right)
\longrightarrow
\frac{1}{2}
\left(
\begin{array}{c}
(M_{11}^2+M_{21}^2)-(M_{12}^2+M_{22}^2) \\[.03in]
2(M_{11}M_{12}+M_{21}M_{22}) \\[.03in]
2|M_{11}M_{22}-M_{21}M_{12}|
\end{array}
\right) 
=
\frac{1}{2}
\left(
\begin{array}{c}
(M^T\! M)_{11}-(M^T\! M)_{22}\\[.03in]
2(M^T\! M)_{12} \\[.03in]
2|\det{M}|
\end{array}
\right) 
=
\frac{1}{2}
\left(
\begin{array}{c}
\cos \lambda  \cos \phi \\[.03in]
\cos \lambda \sin \phi \\[.03in]
\sin \lambda
\end{array}
\right) 
.
$$
This formula may be derived directly or through the Hopf Map described in Section 2.5.5.

A notebook in the programming language Julia is available which takes six representations
and performs all possible conversions.  The notebook provides a correctness check by 
performing all possible roundtrips between representations. (See http://www-math.mit.edu/\verb+~+edelman/Edelman/publications.htm.)

\subsubsection{Disk $r,\phi$  to triangle $a^2,b^2,c^2$ directly}
\label{sub_direct}

An immediate consequence of
$$\left(\begin{array}{c}
a^{2}\\
b^{2}\\
c^{2}\end{array}\right)
=
\mbox{diag}((M\Delta)^T(M\Delta))
= \mbox{diag}(\Delta^{T}V\Sigma^2V^T\Delta)
$$ 
and
$$ 
\Sigma^2=
\left(\begin{array}{cc}
\sigma_1^2 & \\
& \sigma_2^2
\end{array} \right)
=\frac{1}{2} I+r \left(\begin{array}{cc}
1 &\\
& -1 \end{array}\right)
$$ is that 
$$\left(\begin{array}{c}
a^{2}\\
b^{2}\\
c^{2}\end{array}\right)=
\frac{1}{3}
\left(\begin{array}{c}
1\\
1\\
1\end{array}\right)
+r    \left(  \rule{0pt}{0.3in} \Delta_i^T V\left(\begin{array}{cc}
1 &\\
& -1 \end{array}\right)V^T  \Delta_{i} \right)_{i=1,2,3}
.$$
This contains the quadratic form
$V\left(\begin{array}{cc}
1 &\\
& -1 \end{array}\right)V^T$ evaluated at the three columns of $\Delta.$

One can readily  derive 
$$\left(\begin{array}{c}
a^{2}\\
b^{2}\\
c^{2}\end{array}\right)=
\frac{1}{3}
\left(\begin{array}{c}
1\\
1\\
1\end{array}\right)
-\frac{2}{3} r    \left(  \rule{0pt}{0.3in} \begin{array}{l}
\cos \phi_+ \\
\cos \phi_- \\
\cos \phi \end{array}\right)  
\mbox{ with }
\left(  \rule{0pt}{0.3in} \begin{array}{l}
 \phi_+ \\
 \phi_- \\
 \phi \end{array}\right)  
=
\left(  \rule{0pt}{0.3in} \begin{array}{l}
 \phi+{2\pi}/{3} \\
 \phi- {2\pi}/{3} \\
 \phi \end{array}\right)  .
$$

We  prefer the geometric realization that $\Delta$ contains
three columns oriented at angles $\pi/6,5\pi/6,9\pi/6$ with  lengths $\sqrt{2/3}$.
The matrix $V\left(\begin{array}{cc}
1 &\\
& -1 \end{array}\right)V^T$ rotates a vector at angle $\alpha$ clockwise by $\theta$, reflects
across the x-axis, and then rotates counterclockwise by $\theta.$  The quadratic
form thus takes a dot product of a vector at angle $\alpha-\theta$ with its reflection,
i.e., two vectors at angle $2(\alpha-\theta)$, yielding $\cos (2\alpha-2\phi)=-\cos(2\phi-2\alpha+\pi)$.  
Plugging in the three values $\pi/6,5\pi/6,$ and $9\pi/6$ gives the result with little algebra!

As $\phi$ runs from $0$ to $2\pi$, the points $(\sin \phi,\cos \phi)$ trace out a circle in the plane.
The transformation to the plane $x+y+z=0$ is
$$\sqrt{\frac{3}{2}}\Delta^T
\left(   \begin{array}{l}
\sin \phi \\
\cos \phi \end{array}\right)  
=
-
\left(  \rule{0pt}{0.3in} \begin{array}{l}
\cos \phi_+ \\
\cos \phi_- \\
\cos \phi \end{array}\right).  
$$
This yields  an alternative formula to transform  the disk to the squared sidelengths:
$$
\left(\begin{array}{c}
a^{2}\\
b^{2}\\
c^{2}\end{array}\right)=
\frac{1}{3}
\left(\begin{array}{c}
1\\
1\\
1\end{array}\right)
+
\sqrt{\frac{2}{3}}r\Delta^T
\left(   \begin{array}{l}
\sin \phi \\
\cos \phi \end{array}\right)  .
$$

\subsubsection{Triangle $a^2,b^2,c^2$ to disk $r,\phi$    directly: a ``barycentric'' interpretation}
\label{sub_directinverse}

Apply $\Delta$ to the  equation just above:
\begin{equation}
\label{circleq}
\Delta
\left(\begin{array}{c}
a^{2}\\
b^{2}\\
c^{2}\end{array}\right)=
\sqrt{\frac{2}{3}}r
\left(   \begin{array}{l}
\sin \phi \\
\cos \phi \end{array}\right)  .
\end{equation}
Then $r$ and $\phi$ are the polar coordinates of
$$
\tilde{\Delta}\left(\begin{array}{c}
a^{2}\\
b^{2}\\
c^{2}\end{array}\right)= r\left(
\begin{array}{r}
\cos \phi \\
\sin \phi \end{array}
\right), \
{\mbox{ where }} \
\tilde{\Delta}=
\left(
\begin{array}{rrr}
1/2 & 1/2 & -1 \\
\sqrt{3}/2 & -\sqrt{3}/2 & 0 \end{array}
\right) 
 .
$$

{\bf Barycentric Interpretation:}  
To obtain the point on the disk, take an equilateral triangle of sidelength $1$ (the columns of $\tilde{\Delta}$) and find the point with barycentric coordinates $(a^2,b^2,c^2).$  

{\bf The {}``broken stick'' problem:} \cite{goodman08} asks how
frequently $a,b,c$ can be edge lengths of a triangle, when their sum is the length of the stick.
How often do $a,b,c$ satisfy the triangle inequality? This has an easy answer $1/4$ 
(again!) and a long history. 

We could modify the
question to fix $a^{2}+b^{2}+c^{2}=1$. The answer now consists of
those $(a,b,c)$ inside the disk in Figure 7(a) as a percentage of the big triangle.
The solution to our ``renormalized"  problem is that  the fraction $\pi/\sqrt{27}\approx 0.60$ of
all triples $(a,b,c)$ can be edge lengths of a triangle.

\subsubsection{Triangles with fixed area}

In the four representations the following ``sets of fixed area" are  equivalent:
\begin{itemize}
\item Matrices $M$ with ($\|M\|_F^2=\sum M_{ij}^2=1$ and) determinant $\sqrt{12}K$.
\item Triangles with ($a^2+b^2+c^2=1$ and) constant area $K=\frac{1}{4}\sqrt{1-2(a^4+b^4+c^4)}$ (Hero's formula).
\item The latitude at height $\sqrt{12}K$ on the hemisphere.
\item A circle centered at the origin with radius $r=\frac{1}{2}\sqrt{1-48K^2}$ in the disk representation.
\end{itemize}

To sweep through all triangles with area $K$, take  $r=\frac{1}{2}\sqrt{1-48K^2}$
and let $\phi=2\theta$ run from $0$ to $2\pi/3$.  This avoids the ``rotated'' triangles that
cycle through $a,b,c$.

Here are special triangles with area $K$\, (including approximations for very small $K$):
\begin{itemize}
\item {\bf Right Triangles:}

\begin{centering}
\begin{tabular}{|c|c|c|c|} \hline
$K$ & $r$ & $\phi$ (or $\phi_\pm$) & squared sides \\[0.03in]  \hline
$K\le 1/8$ & $r \ge 1/4$& $\pm$acos$(-\frac{1}{4r})$& 
 $ \frac{1}{2}\ $ and 
$\ \frac{1}{4}\pm \frac{\sqrt{1-64K^2}}{4}=\frac{1}{4} \pm \sqrt{\frac{16r^2-1}{3}} $\\ \hline
\end{tabular}
\end{centering}

For area $K\le\frac{1}{8}$ (corresponding to $r\ge\frac{1}{4})$ there are six congruent right triangles.

For $K$ small, the right triangles have squared side lengths $\frac{1}{2}$, $8K^{2}+128K^{4}+O(K^{6})$, and $\frac{1}{2}-8K^{2}-128K^{4}+O(K^{6})$. We can check that $\frac{1}{2}\sqrt{(8K^{2})(1/2)}\approx K$.

\item {\bf Isosceles Triangles} (third side smaller):

\begin{centering}
\begin{tabular}{|c|c|c|c|} \hline
 $K$ & $r$ & $\phi$ (or $\phi_\pm$) & squared sides $\ \ \rule[-.07in]{0in}{0.2in} a^2+b^2+c^2=1$ \\ \hline
$K\le 1/\sqrt{48}$ & $r\le 1/2$ & 0 & 
Two equal sides:$(1+r)/3\ $  Third side:$  (1-2r)/3$
\\ \hline
\end{tabular}
\end{centering}

For each $K$, three congruent isosceles triangles might be called the most acute of the acute triangles (furthest from the three white lines that represent right triangles.) As $r\rightarrow\frac{1}{2},$ this isosceles
triangle approaches a right triangle. (Many physics computations use
this fact.) The squared side lengths for $K$ small are two of size $\frac{1}{2}-4K^{2}-48K^{4}+O(K^{6})$
and a tiny side $8K^{2}+96K^{4}+O(K^{6})$. 

\item {\bf Isosceles Triangles} (third side larger): 

\begin{centering}
\begin{tabular}{|c|c|c|c|} \hline
$K$ & $r$ & $\phi$ (or $\phi_\pm$) & squared sides  $\ \ \rule[-.07in]{0in}{0.2in} a^2+b^2+c^2=1$  \\ \hline
 $K\le 1/\sqrt{48}$ & $r\le 1/2$ & $\pi/3$ & 
Two equal:$(1-r)/3\ $  Third:$  (1+2r)/3$
\\ \hline
\end{tabular}
\end{centering}

For each $K$, three congruent isosceles triangles might be called the most obtuse of the triangles (although they may  be acute).  The squared side lengths for $K$ small are two of size $\frac{1}{6}+4K^{2}+48K^{4}+O(K^{6})$
and a third  side $\frac{2}{3}-8K^{2}+96K^{4}+O(K^{6})$. 

\item {\bf Singular Triangles:} 

\begin{centering}
\begin{tabular}{|c|c|c|c|} \hline
 $K$ & $r$ & $\phi$ (or $\phi_\pm$) & {\bf actual} sides \ \ $a,b,c$ \\ \hline
 $0$ & $1/2$ & any & 
$\sqrt{\frac{2}{3}}|\sin \frac{1}{2}\phi|, \sqrt{\frac{2}{3}}|\sin \frac{1}{2}\phi_+|, \sqrt{\frac{2}{3}}|\sin \frac{1}{2}\phi_-|$
\\ \hline
\end{tabular}
\end{centering}

The longer side is the sum of the two shorter sides.

\item{{\bf Nearly Singular Triangles:}}

For $K$ tiny, $r=\frac{1}{2}-12K^2-144K^4 +O(K^6).$ 

\end{itemize}

We conclude with  area/angle formulas that will be useful in Section 2.5.3.
They may have interest in their own right as they do not seem to appear in the famous Baker
listings of area formulas from 1885. 
\cite{baker85}.

\begin{lem}
The following formulas relate side lengths $a,b,c$, area $K$, and  angles $A,B,C$:
$$
\begin{array}{l}
\tan\, A=4K/(b^{2}+c^{2}-a^{2}) \\
\tan\, B=4K/(a^{2}+c^{2}-b^{2}) \\
\tan\, C=4K/(a^{2}+b^{2}-c^{2}).
\end{array}
$$
\end{lem}
We can derive $\tan A$ from  the law of cosines: $\cos A=(b^{2}+c^{2}-a^{2})/2bc$
and the well known formula $\sin\,A=2K/bc$. We apply Lemma 1 with $a^2+b^2+c^2=1$ in the form
$$
\begin{array}{l}
\tan\, A=4K/(1-2a^{2}) \\
\tan\, B=4K/(1-2b^{2}) \\
\tan\, C=4K/(1-2c^{2}).
\end{array}
$$
 If angle $A$ is held fixed,  there is a linear relationship between $a^2$ and
$K$.  This means that a circular arc on the hemisphere represents all triangles with angle $A$.
This will be further discussed in Section 2.5.3.
Three special arcs (the three white lines in Figure 3) represent all right triangles with
$A$ or $B$ or $C$ equal to $\pi/2$.

\subsection{The Triangle Inequality and the Disk Boundary}

The razor sharp triangle inequality shows up as the smooth disk in the plane $x+y+z=1$.
We saw this in Section 1.2 experimentally.  Careful readers may have thought about this disk
throughout Section 2.3, and some points are worth  attention:

\begin{itemize}

\item {\bf Area Viewpoint:} 
The triangle inequality is equivalent to  non-negative area.
Hero's formula \newline $K=\frac{1}{4}\sqrt{1-2(a^4+b^4+c^4)}$ implies
$a^4+b^4+c^4 \le 1/2$.  This is  the intersection of the sphere of radius $1/\sqrt{2}$ centered at the origin
with the plane $x+y+z=1$ whose distance from the origin measured from $(1/3,1/3,1/3)$ is $1/\sqrt{3}$.
The radius of the circle of intersection is then $\sqrt{1/2-1/3}=\sqrt{1/6}$.
Equation \eqref{circleq} takes us to the disk with $r=1/2$ through the factor of $\sqrt{3/2}$.

For general $a,b,c$ (no restriction on $a^2+b^2+c^2$), Hero's formula is
$$16K^2=(a+b+c)(-a+b+c)(a-b+c)(a+b-c)=(a^2+b^2+c^2)^2-2(a^4+b^4+c^4).$$
This formula has three factors whose positivity amounts to the three triangle inequalities.
The inequality $(x+y+z)^2-2(x^2+y^2+z^2)\ge0$  produces a cone through the origin which intersects
the plane $x+y+z=1$ at the aforementioned circle of radius $1/\sqrt{6}$.

\item {\bf Matrix Viewpoint:}  The boundary of the triangle inequality is equivalent to
$|\det M|=\sigma_1\sigma_2=0.$ 

\item{ \bf Hemisphere Viewpoint:}    det$(M)=0$ corresponds to the
equator of the hemisphere ($\lambda=0$).

\item{\bf Disk Viewpoint:} The conversion table yields  $r=\cos(\mbox{$\lambda$})/2 \le 1/2$, where $\lambda=0$ corresponds to the disk boundary at $r=1/2$.

\item{\bf Triangle Viewpoint:} The sides of the singular triangles (with $a^2+b^2+c^2=1$) are expressed  in the table before Lemma 1.

\end{itemize}

\subsection{Geometric Construction}

Here we exhibit what is marvelous about mathematics\,: \textit{The triangle is not only represented as a point on the hemisphere, but it can be constructed geometrically inside that hemisphere.} 
The full construction will be described in a moment.  Informally, the point on the hemisphere
will be one vertex, and the base  is inscribed in an equilateral triangle on the equatorial plane.

\subsubsection{Parallelians}
\begin{figure}[H]
\begin{centering}
\includegraphics[scale=0.45]{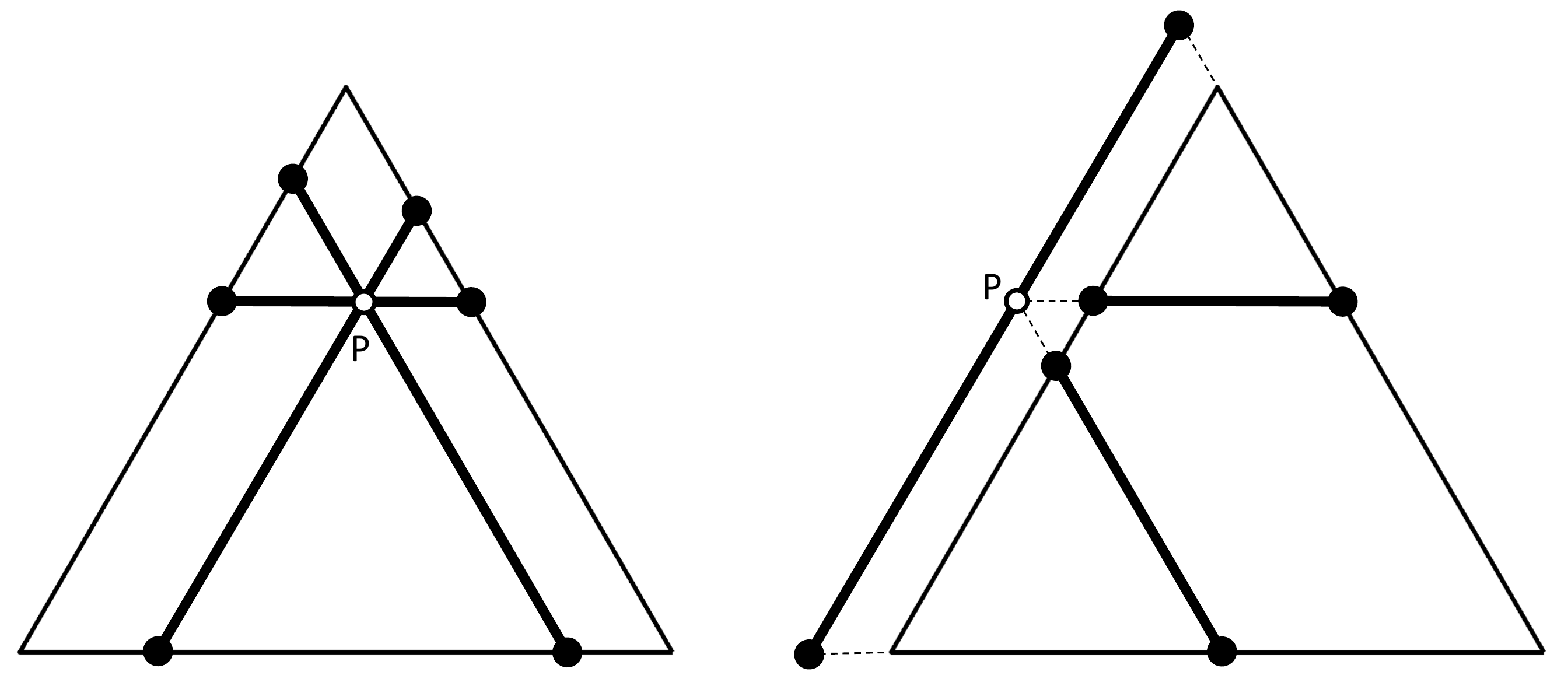}
\par\end{centering}
\caption[]{Parallelians when $P$ is interior or exterior to $E$.}
\label{fig:Parallelians-through}
\end{figure}

We must  define the three parallelians of an equilateral triangle $E$,
meeting at an arbitrary point $P$.
\begin{defn}
\textbf{(Parallelians) }\textsl{Through any point $P$, draw lines parallel to the sides of  $E$ (Figure~\ref{fig:Parallelians-through}). The parallelians are the segments of these three lines with endpoints on the sides (possibly extended) of $E$.}
\end{defn}
When $P$ is interior to $E$, the three line segments are interior
as well. They intersect at $P$. Otherwise their extensions intersect
at $P$.
%
%


\subsubsection{Barycentric coordinates revisited}

This section picks up where Section 2.3.2 left off, with a   
barycentric view of the point in the disk.
We introduce big and little equilateral triangles whose vertices are the columns of
$$\Delta_{\mbox{big}}=
\left(
\begin{array}{rrr}
\sqrt{3}/2 & -\sqrt{3}/2 & 0 \\
1/2 & 1/2 & -1 \end{array}
\right)=\frac{\sqrt 6}{2}\Delta 
$$
$$\Delta_{\mbox{little}}=
\Delta_{\mbox{big}}\left(
\begin{array}{ccc}
0 & 1/2 & 1/2 \\
1/2 & 0 & 1/2 \\
1/2 & 1/2 & 0  \end{array}
\right) =
\left(
\begin{array}{rrr}
-\sqrt{3}/4 & \sqrt{3}/4 & 0 \\
-1/4 & -1/4 & 1/2 \end{array}
\right)
=
-\frac{1}{2}
\Delta_{\mbox{big}}
$$
$$
{\mbox{\hspace{-.7in}or equivalently  }}
\
 \Delta_{\mbox{big}}= \Delta_{\mbox{little}}
 \left(
\begin{array}{rrr}
-1 & 1 & 1 \\
1 & -1 & 1 \\
1 & 1 & -1  \end{array}
\right) .
$$
The columns of $\Delta_{\mbox{big}}$ all have norm $1$.  The midpoints (columns of $\Delta_{\mbox{little}}$)
all have norm $1/2$.

We have  seen in Section 2.4  that if the sides have $a^2+b^2+c^2=1$, then
$$
\Delta_{\mbox{big}} \left(
\begin{array}{r}
a^2  \\
b^2  \\
c^2  \end{array}
\right) =
\Delta_{\mbox{little}} \left(
\begin{array}{r}
-a^2+b^2+c^2  \\
a^2-b^2+c^2  \\
a^2+b^2-c^2  \end{array}
\right) =
\Delta_{\mbox{little}} \left(
\begin{array}{r}
1-2a^2  \\
1-2b^2  \\
1-2c^2 \end{array} \right) 
$$
falls inside a disk of radius 1/2.
(Note when comparing Section 2.3.2, for convenience, we have reflected the x and y coordinates
in $\tilde{\Delta}$.  This serves the purpose of allowing the little triangle to have a horizontal
base.)

In summary, the barycentric coordinates $(a^2,b^2,c^2)$ for the ``big'' equilateral triangle
 correspond to the (possibly nonpositive) barycentric
coordinates $(1-2a^2,1-2b^2,1-2c^2)$ on the inverted ``little'' triangle. (Figure 7(a))

Figure 7  also includes the segments into which parallelians
are sliced.   In an equilateral triangle, the barycentric coordinates are in the same ratio as the
distances to the sides in (d), which in turn are in the same ratio as the parallelian segments  in (e). We note that the circle has radius $1/2$, the little triangle has edges $\sqrt{3}/2$ and the big triangle has edges $\sqrt{3}$.  In Figure 7(e), the lengths labeled  $\alpha,\beta,\gamma$ are exactly $\frac{\sqrt{3}}{2}(1-2a^2),\frac{\sqrt{3}}{2}(1-2b^2),\frac{\sqrt{3}}{2}(1-2c^2)$.

\begin{figure}[H]
\begin{centering}
\includegraphics[scale=0.6]{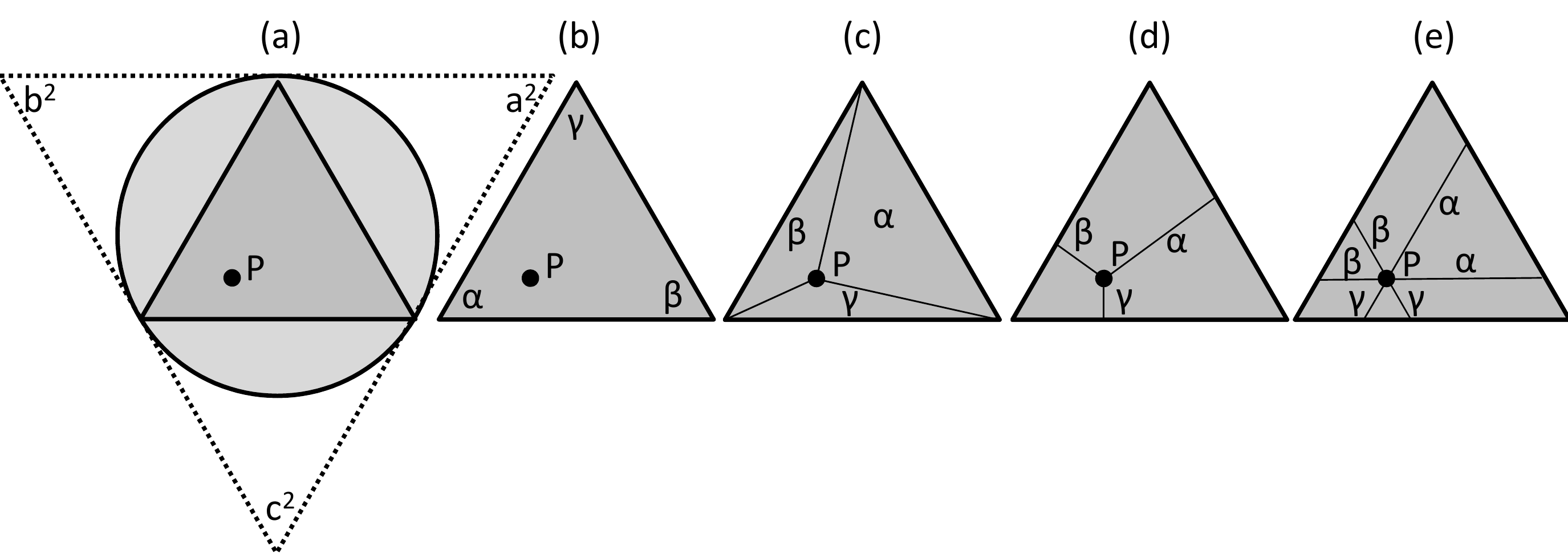}
\par\end{centering}
\caption[]{
(a) $P=(a^2,b^2,c^2)$ in the $\Delta_{\mbox {big}}$ barycentric system.  
(b) $P=(\alpha,\beta,\gamma)$ in the $\Delta_{\mbox{little}}$ system. 
Also in the same ratio $\alpha:\beta:\gamma$ are the areas (c), the perpendicular distances (d), and the parallelian segments (e).
Explicitly $\alpha:\beta:\gamma=(1-2a^2):(1-2b^2):(1-2c^2).$
}
\label{fig:Barycentric-Coordinates-new-in}
\end{figure}

\begin{center}
Endpoints for Parallelians through $P$ in Barycentric Coordinates 
\end{center}

\begin{center}
\begin{tabular}{|c|c|c|} \hline
& $\Delta_{\mbox{big}}$ & $\Delta_{\mbox{little}}$ \\ \hline
$P$ & $(\rule{0pt}{0.2in}a^2,b^2,c^2)$ & $(1-2a^2,1-2b^2,1-2c^2)$ \\[0.1in] \hline
Parallelian 1 &  $\begin{array}{l@{,\ }l@{,\ }l} (\rule{0pt}{0.2in} a^2&1/2&1/2-a^2)\\ (a^2&1/2-a^2&1/2) \end{array}$ 
              &  $\begin{array}{l@{,\ }c@{,\ }r} (\rule{0pt}{0.2in}1-2a^2&0&2a^2)\\ (1-2a^2&2a^2&0) \end{array}$ 
 \\[0.2in]  \hline
Parallelian 2 & $\begin{array}{l@{,\ }l@{,\ }l} (\rule{0pt}{0.2in} 1/2&b^2&1/2-b^2)\\ (1/2-b^2&b^2&1/2) \end{array}$ 
              & $\begin{array}{l@{,\ }c@{,\ }r} (\rule{0pt}{0.2in}0&1-2b^2&2b^2)\\ (2b^2&1-2b^2&0) \end{array}$ 
\\[0.2in] \hline
Parallelian 3 &  $\begin{array}{l@{,\ }l@{,\ }l} (\rule{0pt}{0.2in} 1/2&1/2-c^2&c^2)\\ (1/2-c^2&1/2&c^2) \end{array}$ 
              &  $\begin{array}{l@{,\ }c@{,\ }r} (\rule{0pt}{0.2in}0&2c^2&1-2c^2)\\ (2c^2&0&1-2c^2) \end{array}$
  \\[0.2in] \hline  
\end{tabular}
\end{center}

\subsubsection{Hemisphere to Triangle: Direct Construction}

Let $S$ be a point on the hemisphere corresponding to edge lengths $a,b,c.$
A straightforward but unusual construction yields
a triangle within that hemisphere  whose sides are proportional to $a,b,c$.
In the construction  below, we let $E$ denote the equilateral triangle inscribed
in the hemisphere's equator as in Figure 7(a).  The vertices of $E$ in the xy plane
are the columns of $\Delta_{\mbox{little}}$.
 
\begin{figure}[H]
\begin{centering}
\includegraphics[scale=0.9]{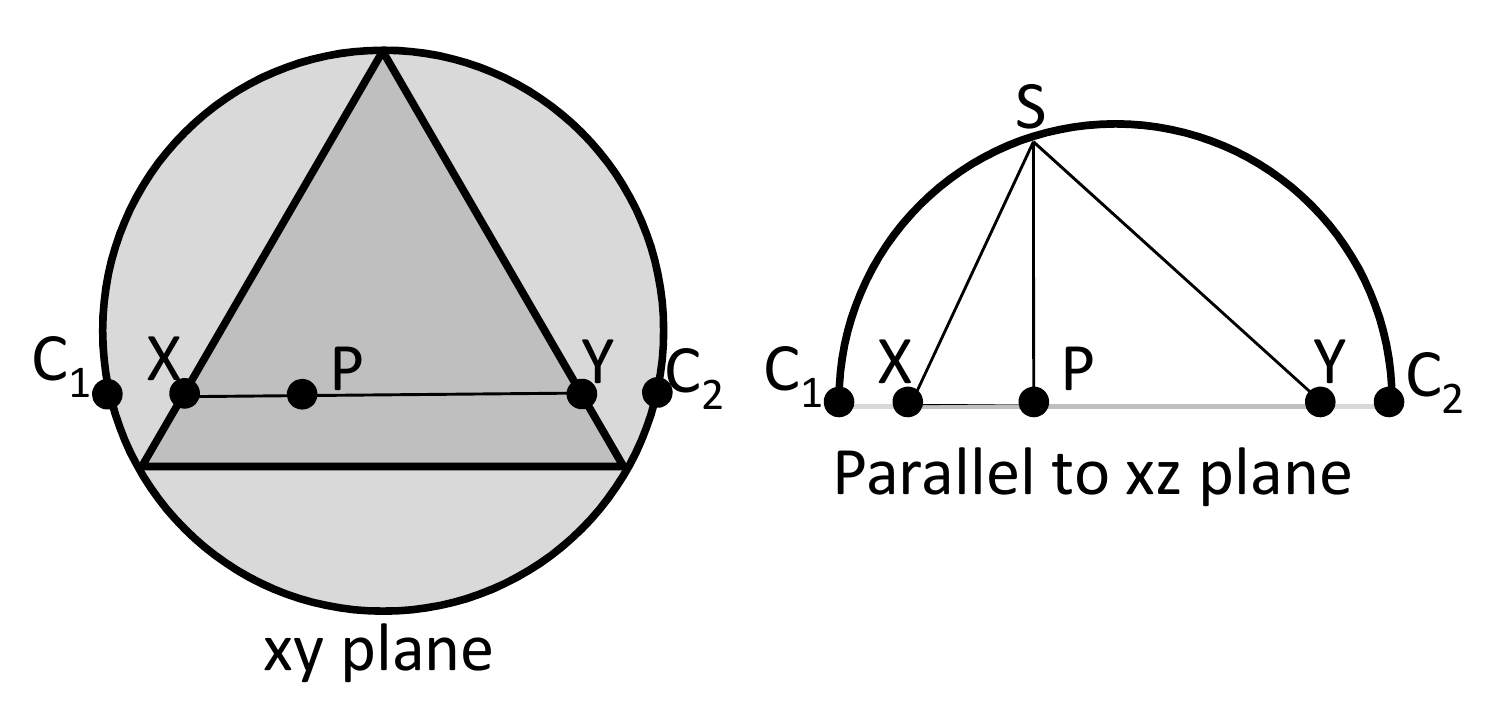}
\par\end{centering}
\caption[]{\begin{tabular}{l} \\
Top View (Left) and Front View (Right) \\
Triangle SXY has sides proportional to $a,b,c$. 
\end{tabular} 
}
\label{fig:construct}
\end{figure}

{\bf Geometry Construction Theorem:} Let S be any point on a hemisphere of radius $1/2$. Let P be its projection onto  the base. The line through P parallel to the x-axis, intersects
$E$ at X and Y.  The triangle SXY in Figure \ref{fig:construct}  has side lengths proportional to $a,b,c$. 

{\bf Proof:} Let $\omega=\sqrt{3}$.  We know from the last sentence in  Section 2.5.2 that XP=$\omega(\frac{1}{2}-b^2)$ and YP=$\omega(\frac{1}{2}-a^2)$. Thus the parallelian length is  XY=$\omega c^2$.  The height SP is $\sqrt{12}K=2\omega K$ from the 
table in Section 2.3.
Thus the tangent of angle SXP is $4K/(1-2b^2)$ and similarly the tangent of angle SYP is $4K/(1-2a^2)$.
These agree with the tangents in Lemma 1, so we have constructed the 
triangle with sidelengths proportional to
$a,b,c.$

The lengths  are   SY=$\omega ac$, SX=$\omega bc$ and XY=$\omega c^2$.
This triangle inside the hemisphere is one of the three similar triangles that we now construct
using parallelians.

The collinear points $P,X,Y$ in the xy plane are the columns of 
$$ [P \ X \  Y]=\Delta_{\mbox{big}} \left(\begin{array}{ccc}
a^2 &  1/2     &  1/2-c^2 \\
b^2 &  1/2-c^2 & 1/2 \\
c^2 &  c^2   & c^2
\end{array}\right) . $$

From
$a,b,c$ we computed
$P,X,Y$ in the equatorial plane in Figure 9 from this formula.
We construct the hemisphere point $S$  by changing $P$'s
 z-coordinate to  $\sqrt{\frac{1}{4}-\|P\|^2}$.

\subsubsection{The three similar triangles theorem}

Theorem 3 is interesting for several reasons.  It seems unusual
that in traditional geometry, one would construct triangles with vertex on a hemisphere, and opposite edges inscribed in an equilateral triangle that itself is inscribed in the hemisphere's base.
It also feels unusual that the three triangles would share a common line segment SP which represents
{\it three different altitudes} in the shape.

\begin{thm}
S is on the hemisphere corresponding to $a,b,c$  and P is its projection onto the  base.  Three triangles
share SP as an altitude with vertex S and opposite edge one of the parallelians 
through P as in Figure 7(e). The three triangles thus created are similar with  edges in the proportion $a:b:c$. 
\end{thm}

{\bf Proof:}  The construction in 2.4.3, by symmetry, yields three  triangles. The scaling factors $2\omega a$ and
$ 2\omega b$ and  $2\omega c$ multiply $a,b,c$.
The three similar triangles are illustrated in Figure 9.

\begin{figure}[H]
\begin{centering}
\includegraphics{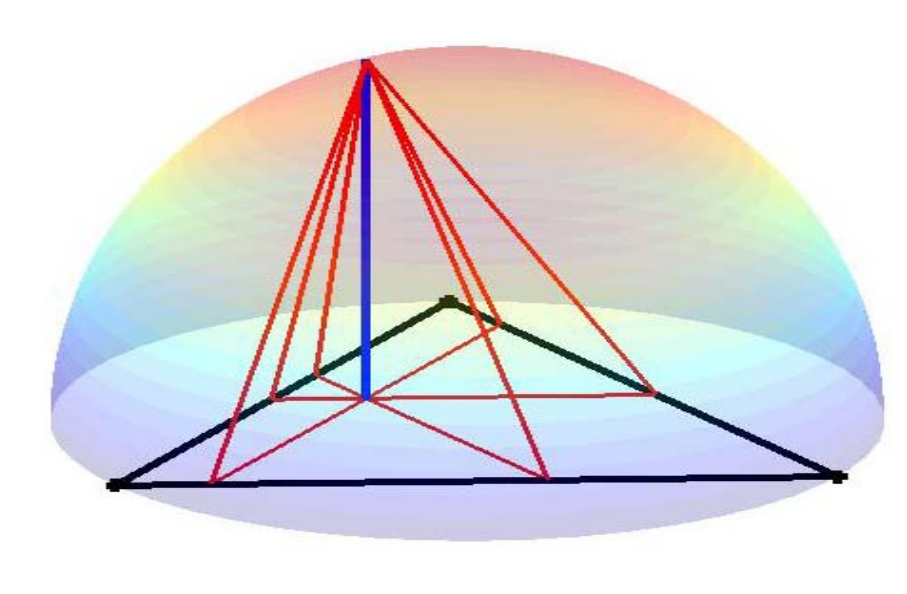}
\caption[]{
\begin{tabular}{l} \\ 
Illustration of the hidden geometry inside triangle shape space revealed in Theorem 3: \\
The three similar triangles theorem. \\ 
Above: Skeleton: The three triangles sharing the altitude SP are similar\,! \\
Below: 3d model: All triangles are similar; 
congruent triangles are color coded.
\end{tabular}
}
\vspace{.2in}
\begin{center}
\includegraphics[scale=.35]{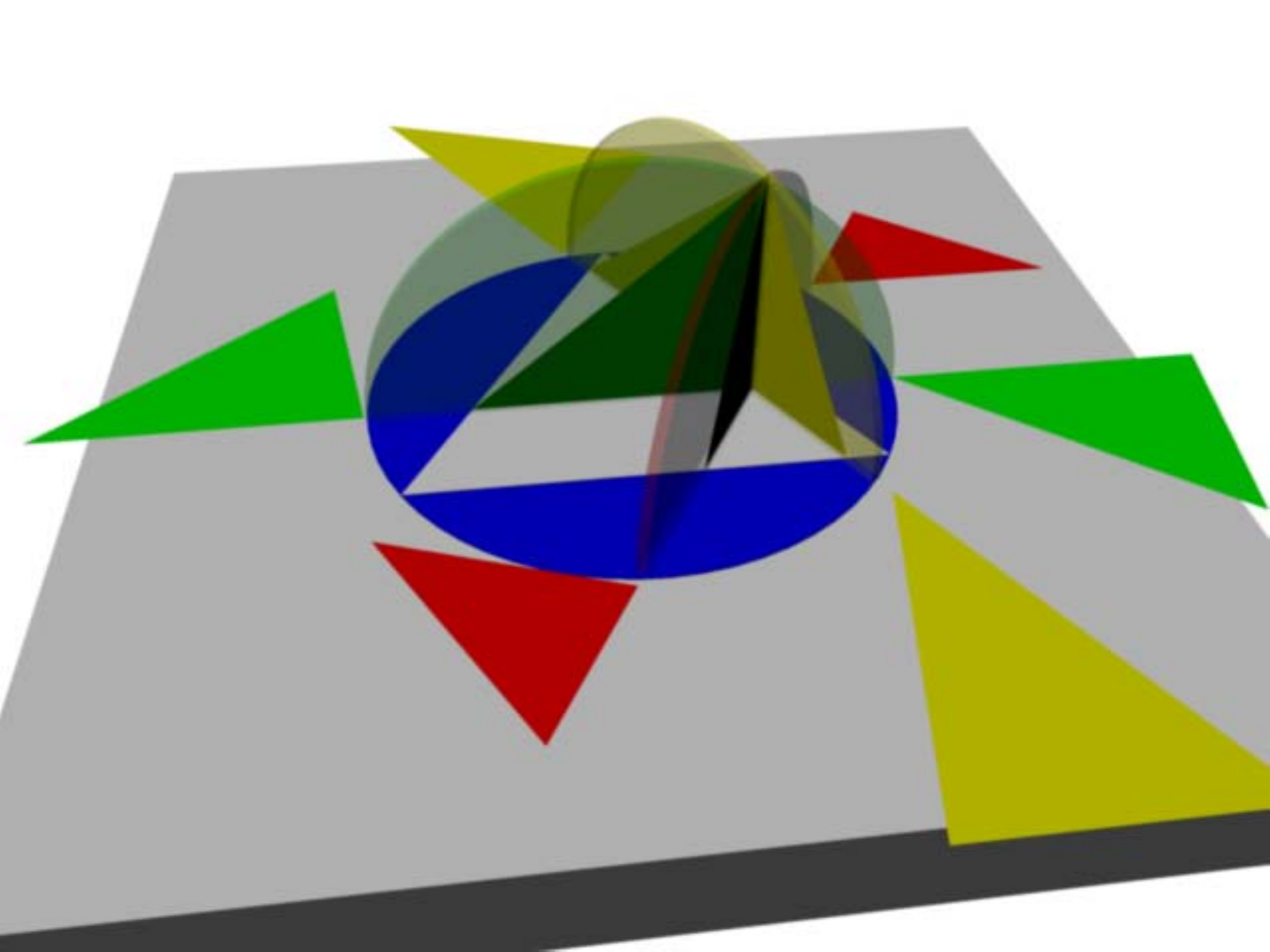}
\includegraphics[scale=.35]{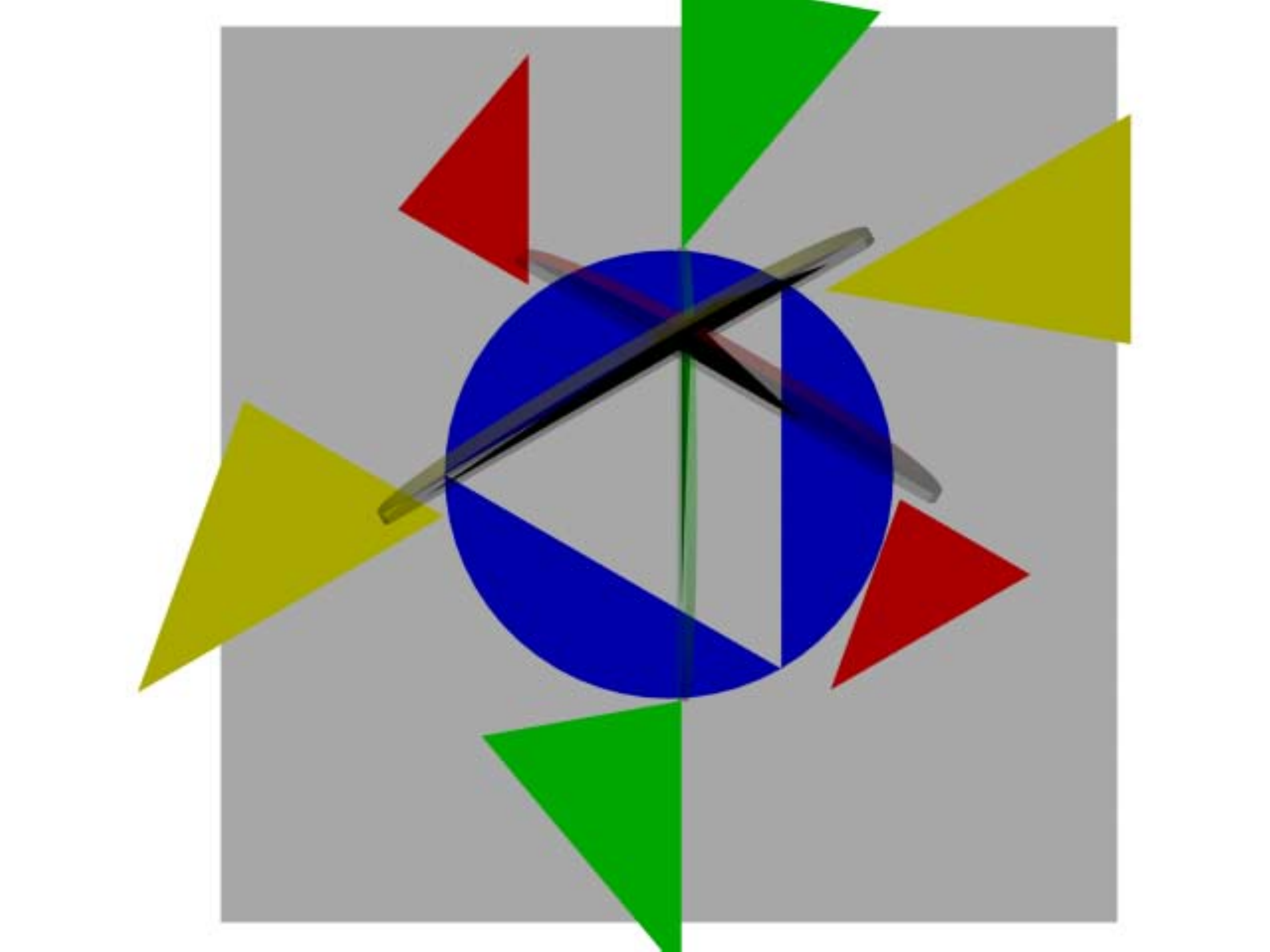}
\end{center}
\end{centering}
\end{figure}

\newpage

\subsubsection{Matrix to Hemisphere: The Hopf Map}

Hopf discovered in 1931  (see \cite{Lyons_2003}) that the hypersphere $S^3$ maps naturally to the sphere $S^2$.
We think of points on $S^3$ as the $2\times2$ matrices normalized by $\|M\|_F^2=\sum M_{ij}^2=1$.
The Hopf map comes directly from the Singular Value Decomposition  $M=U\Sigma V^T$:
$$
\mbox{Hopf($M$)}:
M=U
\left(
\begin{array}{cc}
\cos(\lambda/2) & 0 \\
0 & \sin(\lambda/2)
\end{array}
\right)
\left(
\begin{array}{cr}
\cos(\phi/2) & -\sin(\phi/2) \\
\sin(\phi/2) & \cos(\phi/2)
\end{array}
\right)^T
\longrightarrow
\left(
\begin{array}{l}
\cos \lambda \cos \phi \\[.03in]
\cos \lambda \sin \phi \\[.03in]
\sin \lambda
\end{array}
\right) .
$$
in $S^2$.

Here the latitude is  $\lambda \in [-\pi/2,\pi/2]$, with sign matching $\det(M)$.
The  $|\lambda|$ is determined by the singular values of $M$, while the longitude $\phi$ comes from the
right singular vectors of $M$.  The rotation matrix $U$ of left singular vectors gives the fiber of $S^3$, that is mapped
to a single point on $S^2$.

Remark: We have not seen the SVD in the usual definitions of the Hopf fibration.  It may be the the most immediate
way to get to the Hopf map.

An explicit formula is
$$
\mbox{Hopf($M$)}:
\left(
\begin{array}{cc}
M_{11} & M_{12} \\
M_{21} & M_{22}
\end{array}
\right)
\longrightarrow
\left(
\begin{array}{c}
(M_{11}^2+M_{21}^2)-(M_{12}^2+M_{22}^2) \\[.03in]
2(M_{11}M_{12}+M_{21}M_{22}) \\[.03in]
2(M_{11}M_{22}-M_{21}M_{12})
\end{array}
\right) .
$$

Our construction of triangle space, mapping
$M$ to the upper hemisphere,  uses $\frac{1}{2}$Hopf($M$). The  coordinate, $z=\frac{1}{2}\sin \lambda$, is
always chosen positive.

\subsubsection{Rotations, Quaternions, and Hopf}

The bad news is that 
there is an ugly formula for the general 3 $\times$ 3 rotation matrix:
$$ 
Q_3 =
\left( \begin{array}{ccc}
\alpha^2 +\beta^2-\gamma^2-\delta^2 & 2(\alpha\delta+\beta\gamma) & 2(\beta\delta-\alpha\gamma) \\
-2(\alpha\delta-\beta\gamma) & \alpha^2 -\beta^2+\gamma^2-\delta^2 & 2(\alpha\beta+\gamma\delta) \\
2(\alpha\gamma+\beta\delta) & -2(\alpha\beta-\gamma\delta) & \alpha^2-\beta^2-\gamma^2+\delta^2 
\end{array} \right) , $$
where $1=\alpha^2+\beta^2+\gamma^2+\delta^2$.
This is a rotation about 
axis $(\beta,\gamma,\delta)$
with rotation angle $2$  acos($\alpha$). 

The good news is that given any $Q_3$
we can obtain $\alpha,\beta,\gamma,\delta$ 
through an eigendecomposition.  We can then associate a 4$\times$4 rotation matrix
$$Q_4=\left( \begin{array}{rrrr}
\alpha & -\beta & \delta & -\gamma \\
\beta & \alpha & \gamma & \delta \\
-\delta & -\gamma & \alpha & \beta \\
\gamma & -\delta & -\beta & \alpha 
\end{array} \right) .
$$

The key relationship that ties $Q_3$ to $Q_4$  will be crucial for Theorem 6:
$$\mbox{Hopf}(Q_4M) = Q_3 \mbox{Hopf}(M).$$
$Q_4M$ denotes the operation of flattening $M$ to the column vector
$[M_{11},M_{21},M_{12},M_{22}]^T,$ applying $Q_4$, and reshaping back
into a two by two matrix.

Readers who prefer not to think about quaternions can safely skip ahead to Section 2.6.
$\mbox{Hopf}(Q_4M) = Q_3 \mbox{Hopf}(M),$ could be checked directly, for example by Mathematica.
But  quaternions are the better way to understand why the
relationship holds. It
is  a manifestation of the quaternion identity 
$$(qm)i(\bar{qm})=q(mi\bar{m})\bar{q}.$$

We encode the matrix $M$ with the unit quaternion $m=M_{11}-M_{21}i-M_{22}j+M_{12}k$, and
we will encode $Q_3$ and $Q_4$ with the unit quaternion $q=-\alpha+\beta i +\gamma j +\delta k.$
The Hopf map is visible in the relationship:
$$
\begin{array}{rcl}
mi\bar{m}&=&  \left((M_{11}^2+M_{21}^2)-(M_{12}^2+M_{22}^2) \right)\ i \  +  \
 2(M_{11}M_{12}+M_{21}M_{22})\   j \  +\   2(M_{11}M_{22}-M_{21}M_{12})\  k  \\[.07in] 
& = &[i \ j \ k ]\ \mbox{Hopf($M$)}.  \end{array}$$

The matrix $Q_3$ may be found in the computation of $q(iX+jY+kZ)\bar{q}$ whose $(i,j,k)$ parts
are the components of $Q_3[X\  Y\ Z]^T$.

The matrix $Q_4$ may be found in the computation of  the quaternion product  $w=qm$ by writing
$W=
\left(
\begin{array}{rr}
w_1 & w_4 \\
-w_2 & w_3
\end{array}\right).$
This reverses  how $m$ was formed from $M$. We then create $Q_4$ from $W=-Q_4M$ by flattening $M$ and $W$.

\subsection{The Hemisphere: Positions  of special triangles}

We conclude the nonrandom section of this paper with a  summary of the correspondence
between position on the hemisphere and special triangles:

\begin{itemize}
\item{Lines of Latitude:} Triangles of equal area.
\item{Lines of Longitude:} Triangles resulting from a fixed rotation of the reference triangle followed by  scaling along the x and y axes.
\item Circular arc in a plane through an edge of the little equilateral triangle $E$: Triangles with constant angle. (Around the perimeter of the circle, most points correspond to the degenerate triangles with angles $(0,0,\pi)$.  At the vertices of $E$, we obtain all the triangles $(\theta,\pi-\theta,0)$.
\item{Vertical circular arcs as above:}  Right triangles.
\item{North Pole:}  Equilateral triangle.
\item{Longitudes at multiples of $\pi/3$:} Isosceles triangles
\item{The Equator:} Singular triangles.
\item{Inner Spherical Triangle:} Acute triangles 
\item{The six fold symmetries:}  The six permutations of unequal $a,b,c$. 
\item{Where is the lower hemisphere?:} Replace $\det M$ with $-\det M$. Every ordered triangle is ``double
covered."  A signed SVD with rotations $U$ and $V$ would work well.
\end{itemize}

\section{The Normal Distribution generates random shapes}

Section 2  concentrated on the relationship between triangles,
the hemisphere and disk, and 2x2 matrices through the SVD.  
Nothing was random.  Here in Section 3,
randomness comes into its own with the normal distribution as a natural structure.

\subsection{The Normal Distribution}

The normal distribution has a very special place in mathematics.  
Several decades of research into Random Matrix Theory have  shown that exact analytic formulas
are available for important distributions when the matrix entries arise from independent
standard Gaussians.  Non-Gaussian  entries rarely enjoy the same beautiful properties for finite matrices,
though recent ``universality" theorems show convergence for infinite random matrix theory.

We might anticipate, if only by analogy, that random triangles arising from independent
Gaussians might also enjoy special analytic properties compared to other distributions.
Section 3.2  reveals that random triangles generated
from a Gaussian random matrix are uniformly distributed on the hemisphere.

\vspace{.1in}

{\bf Four Gaussian degrees of freedom  from six:}

Figure 2 illustrates 1000 random triangle shapes. 
As we saw in Section 1.3, six independent standard Gaussians may be used 
to describe the 2x3 array of vertex coordinates with  six degrees of freedom.

Nonetheless, if we translate the  centroid to the origin, or a vertex to the origin, or use the $2 \times 3$ edge matrix $E$,
or the matrix $M_e$ from the edge viewpoint or the matrix $M_v$ from the vertex viewpoint, only 
four degrees of freedom remain.  

Consider the $2 \times 3$ matrices $T$, the origin centered vertex matrix, or $E$ the edge matrix. 
Up to scaling, both of these  matrices consist of  elements that are standard normals conditioned on a zero column sum.
This may be argued or computed.  The argument is that Gaussians are linear, and the symmetry of the situation requires
that no one vertex or edge can be special.

One can check readily  that centering  six standard normals
at the origin  produces two independent rows  each with the singular covariance matrix
$$\Delta^T  \! \Delta=\frac{1}{3}
\left(
\begin{array}{rrr}
2 &-1 &-1 \\
-1& 2 &-1 \\
-1 &-1& 2
\end{array}
\right) .
$$
The matrix $\Delta^T  \! \Delta$ is the projection matrix onto the plane $x+y+z=0$.
The  edge matrix $E$ obtained by taking differences of columns  can be readily verified to have
covariance matrix  $3\Delta^T  \! \Delta. $  As shape information removes  scale information,
this is irrelevant to the paper, and we often tend to think of $T$ and $E$ as providing dual views of the 
same construction.

If one forms the matrix $M_v=T\Delta^T$, the covariance matrix for each row of two elements  is now $ \Delta (\Delta^T \! \Delta) \Delta^T = I_2.$
We then obtain the nice conclusion that $M_v$ consists of four independent standard normals. Similarly,  $M_e$ before scaling would thus
consist of four independent normals with  variance $3$.
Once again we remark that scaling is irrelevant to shape theory results justifying our normalization of choosing $M$ to have sum of squares $1$.

\subsection{Uniform distribution on the hemisphere}

\begin{prop}
These equivalent representations describe
the uniform measure on the hemisphere with radius $1/2$:
\end{prop}

\begin{tabular}{|l|l|} \hline
Longitude & 
$\phi$ uniform on $[0,2\pi)$ is independent from: 
\\ \hline \hline
Height & uniform on $[0,1/2]$         
\\ \hline
 Latitude &
    density $\rule{0in}{0.15in}\cos(\lambda)d\lambda \mbox{ on } [0,\pi/2]$        
\\ \hline
$r$ &  density
 $\rule{0in}{0.15in} 
4r(1-4r^2)^{-1/2}dr \mbox{ on } [0,1/2] $                
 \\ \hline
 $|\det M|$ & uniform on $\rule{0in}{0.15in}[0,1/2]$ \\ \hline
 Triangle area & uniform on $\rule{0in}{0.15in} [0,1/\sqrt{48}]$  \\ \hline
\end{tabular}

\vspace{.15in}

{\bf Proof:}  The latitude and longitude joint density  $\cos(\lambda)d\lambda d\phi$ is the familiar volume element on the sphere,
where latitude $\lambda$ is measured from the equator rather than the zenith.
Archimedes of Syracuse 
(born -287, died -212) knew the height  formula
through his hat box theorem: on a sphere, the surface area of a zone is proportional to its height.  
The area comes from Hero (c.\ 10-70)  or Heron of Alexandria's formula, believed already known to
Archimedes centuries earlier.

\begin{lem}
{\bf Exponential distributions}  The sum of squares of two independent standard normals ($\chi_2^2$)
is exponentially distributed with  $\mbox{density } \frac{1}{2}e^{-x/2}$.  If $e_1,e_2$ are  independent 
random variables with identical exponential distribution then $e_1/(e_1+e_2)$ is uniformly
distributed on $[0,1]$.
\end{lem}

{\bf Proof:}  These facts are well known.  The generalization to $n$ exponentials
is a popular way to generate uniform samples $x_i=e_i/\sum e_i$ on the simplex $\sum_{i=1}^n x_i=1, x_i\ge 0$
(cf. \cite[Proposition 1]{Portnoy94}).  

\vspace{.15in}

We now turn to the beautiful result known to David Kendall and his collaborators:

\begin{thm}
Triangles generated from six independent normals ( xy coordinates of three vertices or edges) correspond
to points distributed uniformly  on the hemisphere.
\end{thm}

The corresponding $M$ is a $2 \times 2$ matrix of independent standard normals.

{\bf Proof 1 (Exponentials):}
We first consider the height:
$$\det M/\|M\|_F^2 =\frac{ad-bc}{a^2+b^2+c^2+d^2}=
\frac{
   \left( \left( \frac{a+d}{\sqrt{2}} \right)^2+
     \left( \frac{b+c}{\sqrt{2}} \right)^2 \right) -
    \left(  \left( \frac{a-d}{\sqrt{2}} \right)^2+
       \left( \frac{b-c}{\sqrt{2}} \right)^2 \right) }
       {2\left( \left( \frac{a+d}{\sqrt{2}} \right)^2+
     \left( \frac{b+c}{\sqrt{2}} \right)^2  +
     \left( \frac{a-d}{\sqrt{2}} \right)^2+
       \left( \frac{b-c}{\sqrt{2}} \right)^2 \right)       
       }
    . $$


The numerator  is the difference
between two exponentially distributed random variables
(being the sums of squares of two independent standard normals.)
Then by Lemma 5,
$$\frac{e_1-e_2}{2(e_1+e_2)}=
\frac{e_1}{e_1+e_2}-\frac{1}{2}$$ is uniform on $[-1/2,1/2]$.
The absolute value is then uniform on $[0,1/2]$ as desired. 

If we now turn to the ``longitude,"
the right singular vectors provide the uniform
distribution  owing to the right orthogonal invariance of $M$.

\vspace{.1in}

{\bf Proof 2 (Random Matrix Theory and condition numbers):}

The singular values $\sigma_1,\sigma_2$ and singular vectors $V$ for  $M$=\verb+randn(2,2)+  are independent.  $V$ is a rotation with angle uniformly distributed on $[0,\pi).$

 The distribution of the condition number $\kappa$ is important. It may be found  in
 \cite[Eq. 2.1]{edelman88a} or \cite[Eq. 14]{edelman89a}.
Restated, the probability density of $\kappa$ for a random $2\times2$ matrix of iid normals is $-2\tfrac{d}{dx}(x+x^{-1})^{-1}$. Equivalently, the distribution of
$
\frac{\sigma_{1}\sigma_{2}}{\sigma_{1}^{2}+\sigma_{2}^{2}}=(\kappa+\kappa^{-1})^{-1}
$
is \textbf{uniform} on $[0,1/2]$.

\vspace{.1in}

{\bf Proof 3 (Hopf Map)}

We know that $M$ is uniformly distributed on the sphere $1=M_{11}^2+M_{21}^2+M_{12}^2+M_{22}^2.$
We seek to explain why Hopf($M$) is uniform on the sphere $1=x^2+y^2+z^2$.
The mathematician's favorite proof would have the latter  inherited from 
the former.  We want  uniformity of Hopf$(M)$  to be a shadow of uniformity of $M$, but the Hopf map is not
a linear projection.

For any 3$\times$3 rotation matrix $Q_3$, we  constructed in Section 2.5.6  a
4$\times$4 rotation matrix $Q_4$, such that  
$$\mbox{Hopf}(Q_4M) = Q_3 \mbox{Hopf}(M),$$
for any fixed matrix $M$.
If $M$ is random and uniformly distributed on
$1=M_{11}^2+M_{21}^2+M_{12}^2+M_{22}^2,$ then of course $Q_4M$
is too.

What does this say about the distribution of Hopf($M$)?
For any choice of $Q_3$, the distribution of $Q_3\mbox{Hopf}(M)$
is the distribution of Hopf($Q_4M$), which is the same as the distribution of Hopf($M$).
The distribution of $\mbox{Hopf}(M)$ is invariant
under any fixed 3$\times$3 rotation.
It must be the uniform distribution on $x^2+y^2+z^2=1$.

We encourage the reader to follow this -- it is truly a beautiful argument.

%
%

\subsection{What is the probability that a random triangle is acute or obtuse\,?}

Put  Lewis Carroll's question in the context of the normal distribution.   We then obtain the result
known to Portnoy and  Kendall and collaborators:

\begin{thm}
A random triangle from the uniform distribution in shape space has squared side lengths uniformly distributed on $[0,2/3]$.
The probability 
is $1/4$ that this triangle is acute.
\end{thm}

{\bf Proof:}
The edges of the triangle are the lengths of the three columns of $M\Delta$.  Taking the third column
for convenience, the squared side length is $\Delta_3^2=2/3$ times $(M_{21}^2+M_{22}^2)$.
By Lemma 5, this is the uniform distribution on $[0,2/3]$.

Edges that satisfy $a^{2}+b^{2}+c^{2}=1$ give a right triangle  when $c^{2}=1/2$. The triangle is obtuse when $c^{2}>1/2$, which makes $c$ the longest side. The probability that a particular angle is obtuse is then $(2/3-1/2)/(2/3)=1/4$. The probability that any angle is obtuse is then $3/4$ (at most one can be obtuse!). Then $1/4$ is the probability that all are acute.\qed

\subsection{Triangles in $n$ dimensions}

There is an obvious generalization of Figure 5 to higher dimensions.  Let $M$ be a random matrix of independent
standard norms with $n$ rows and $2$ columns.  Then $M\Delta$ is a random triangle  shape in $n$ dimensions.
We can think of this as random vertices centered at the origin, or random edges that close to a proper triangle.
There is a further generalization

The same argument   as for the plane shows a squared side length has the distribution $(2/3)\mathrm{Beta}(n/2,n/2)$. In $\mathbb{R}^{n}$    the probability of an obtuse triangle  is
 $3(1-I(\frac{3}{4},\frac{n}{2},\frac{n}{2}))$ where $I$ denotes the incomplete beta function. This can be evaluated in MATLAB by
 $$\mbox{\texttt{3{*}(betainc(3/4,n/2,n/2,'upper'))}}$$ or in Mathematica
by $$\mbox{\texttt{3{*}(1-BetaRegularized[3/4, n/2, n/2])}.}$$

This probability is also computed by Eisenberg and Sullivan 
\cite{eisenberg96}. They note that for larger $n$, it is increasingly likely that a random triangle is acute. The probability of an obtuse triangle is $10\%$ in $\mathbb{R}^{12}$, and $1\%$ in $\mathbb{R}^{26}$. It is nearly $0.1\%$ in $\mathbb{R}^{40}$.

For arbitrary $n$  the distribution is no longer uniform on the hemisphere. 

Connecting multivariate statistics and numerical
analysis, the same result is innocently hidden as an exercise
in Wishart matrix theory.  It is an unlikely connection. This ``square root ellipticity statistic''  may be found in Exercise~$8.7(\mbox{b})$ of 
\cite[p.379]{Muirhead82a}. It states that \foreignlanguage{english}{$P\left(\frac{2\sigma_{1}\sigma_{2}}{\sigma_{1}^{2}+\sigma_{2}^{2}}<x\right)=x^{n-1}$}. As $n$ increases, this reduces the probability near the equator and adds weight near the poles. Acute triangles become more
probable.  

\subsection{An experimental investigation that revealed the distribution's {}``shadow''}
Early in our own investigation, we plotted the normalized squares of $\mathtt{randn(2,2)}*\Delta$ as barycentric coordinates. We found very quickly that the probability of an acute triangle is $1/4$, and triangles naturally fill a hemisphere with a uniform distribution. Here is a MATLAB code and a picture that tells so much in so little space that we could not resist sharing. A few interesting things to note\,: Line $10$ projects the three  normalized squared edges to $\mathbb{R}^{2}$ using $\Delta\transp$. Line 9 compares $a^{2}+b^{2}$ to $c^{2}$ to decide acute or obtuse. (This is more efficient than computing angles.) Ten thousand trials gave $25.37\%$ acute triangles. Ten million trials gave $24.99\%$ acute triangles.

Readers may wish to recreate this road to discovery of the uniform
distribution on the hemisphere\,: Histogram the radii of the points, guess the measure, and then realize the uniformity. We guessed that uniformity by fitting the density $f(x)=4x/\sqrt{1-4x^{2}}=-\frac{d}{dx}\sqrt{1-4x^{2}}$ which is the
shadow of the uniform distribution on the hemisphere (differential form
of Archimedes' hat box theorem). Then came proofs.
\vspace{-.17in}
\begin{figure}[H]
\includegraphics[scale=1.0]{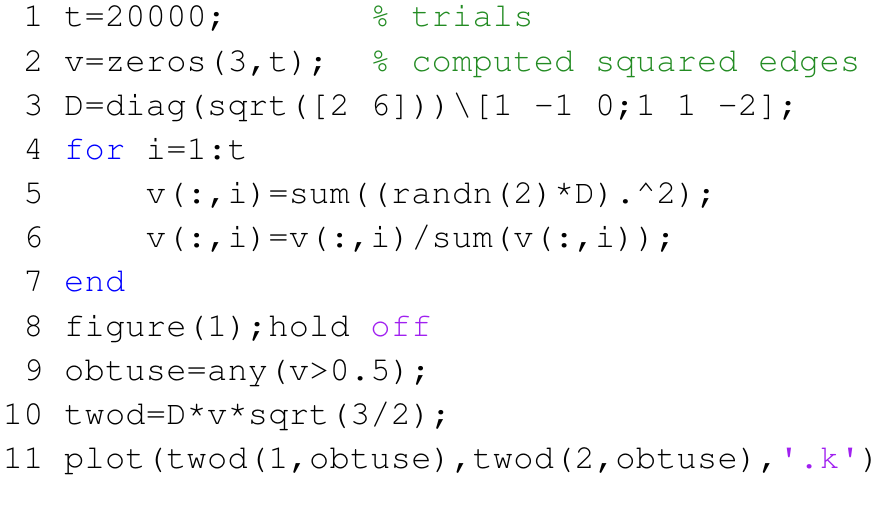}
\includegraphics[scale=.4]{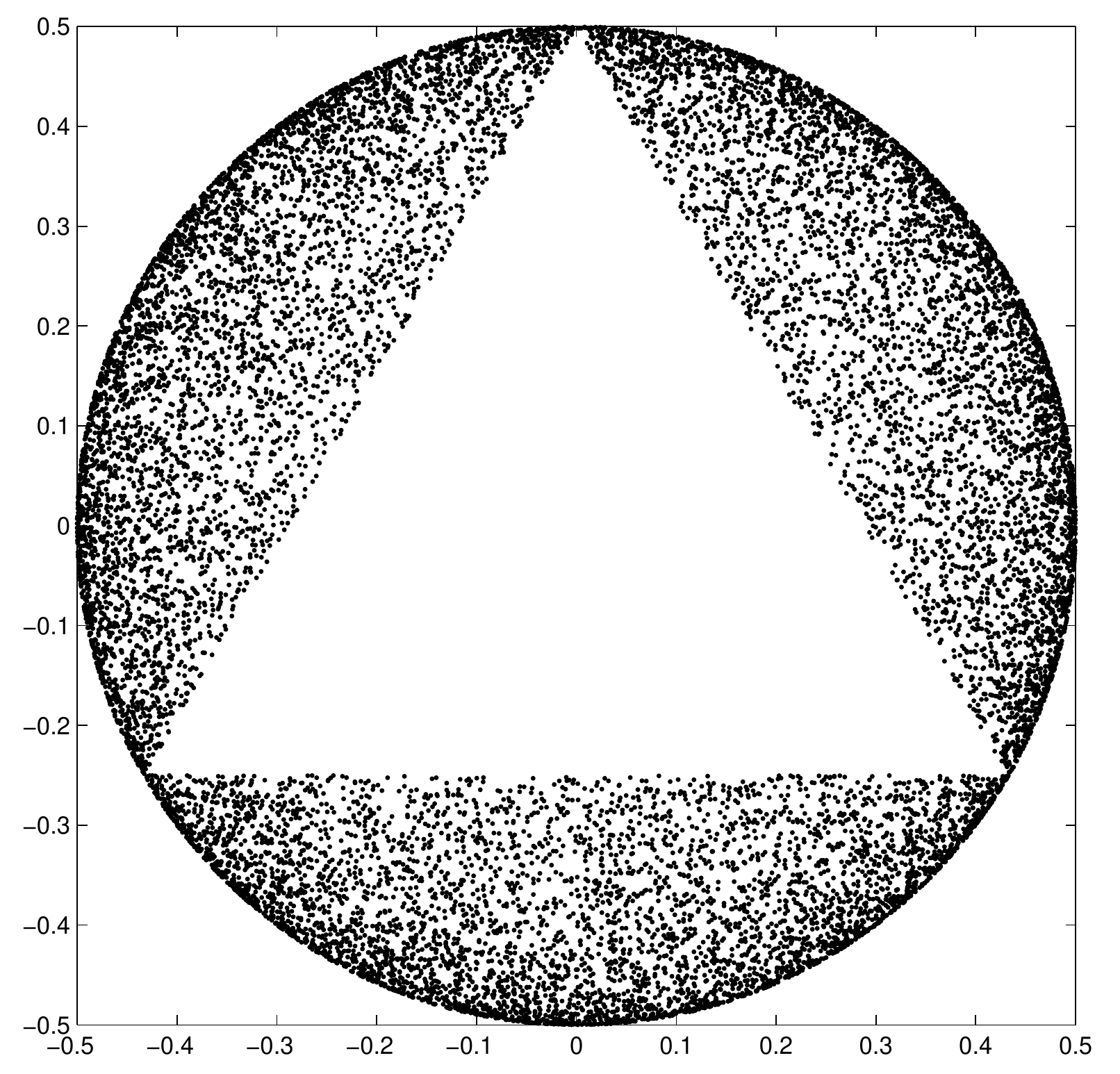}
\caption[]{Random $a^2,b^2,c^2$  are computed on line $5$ and normalized on line $6$. The obtuse triangles are identified on line $9$. Line $10$ projects the plane $x+y+z=1$ onto $\mathbb{R}^2$. Line $11$ plots every obtuse triangle as a point in this plane. The   acute triangles form the inner triangle as in Figure $3$ and  the obtuse triangles  complete a disk of radius $1/2$. The hemisphere is viewed from above.}
\end{figure}

\vspace{-.3in}
\begin{figure}[H]
\includegraphics[scale=1.0]{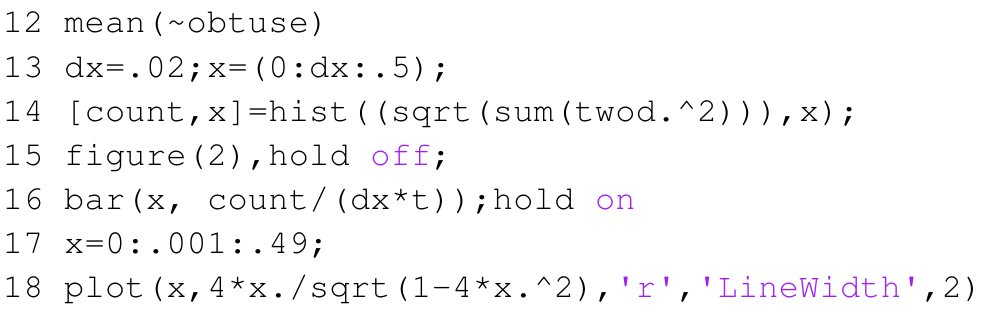}
\includegraphics[scale=0.27]{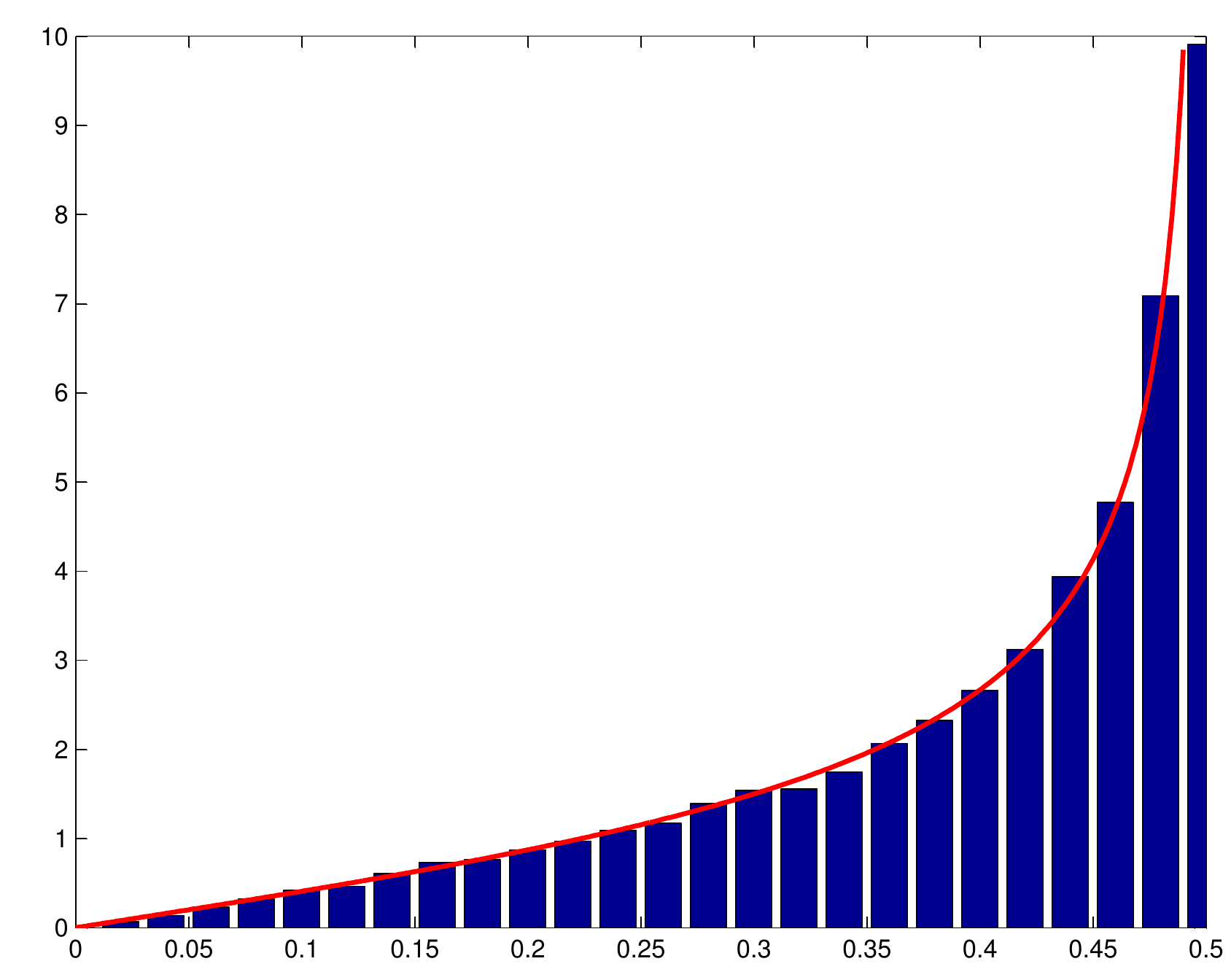}
\vspace{0.2in}
\caption[]{Line $12$ computes the sample probability of acute triangles. Lines $13$ through $18$ histogram all the radii against the guess (curved line) that the points are uniformly distributed on the hemisphere.}
\end{figure}

\subsection{Uniform shapes versus uniform angles}
Two uniform distributions, on the hemisphere and on angle space $A+B+C=180^{\circ}$, gave the same fraction $\frac{3}{4}$ of obtuse triangles. We are not aware of a satisfying theoretical link. Portnoy 
\cite[Section 3]{Portnoy94} philosophizes about {}``the fact that the answer $3/4$ arises so often''. To emphasize the difference between these  distributions, we report on a numerical experiment and a theoretical density computation to understand the angular distribution.

Our first figure might be called $100\!,000$ triangles in $100$ bins. The three angles divided by $\pi$ are barycentric coordinates in the
figure. With four bins the triangles would appear uniform. With $100$ bins we see that they are anything but. A uniform distribution would have 1000 triangles per bin.
\begin{figure}[H]
\begin{center}
$100,\!000$ triangles in $100$ bins
\end{center}

\begin{center}
\includegraphics[width=300pt]{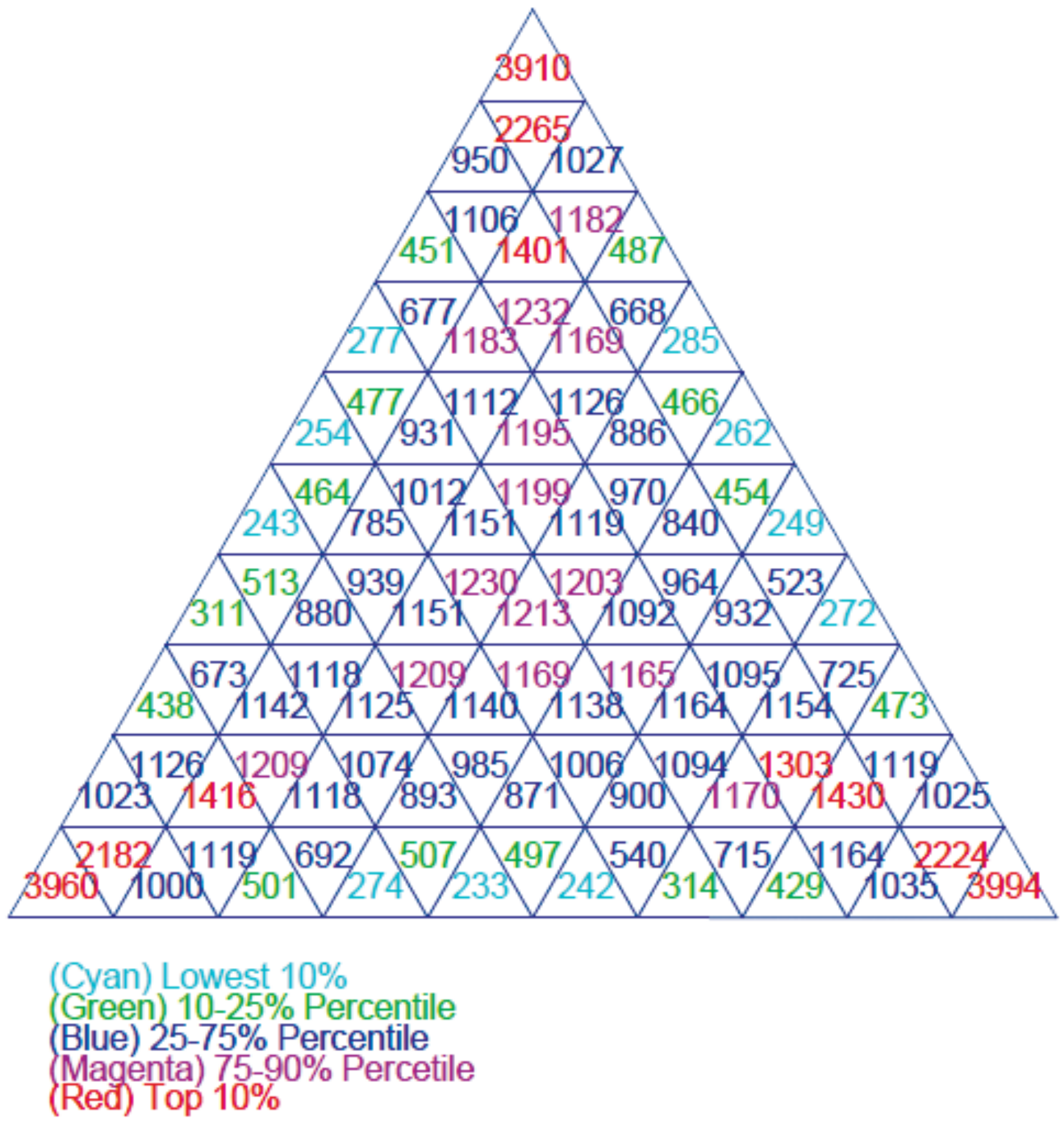}	
\end{center}
\caption[]{$100\!,000$ triangles in $100$ bins\,: Triangles selected uniformly in shape space are not uniform in angle space.}
\end{figure}

The underlying theoretical distribution involves 
the Jacobian from angle space to the ``squared edges'' disk. Then
one must apply a second Jacobian to the hemisphere, and finally
invert. Since the second part is standard, we will sketch the first
part of the argument. 

Suppose  $\alpha+\beta+\gamma=1$. Considering
the Law of Sines, define
%
\[
s=\left(\begin{array}{c}
a^{2}\\
b^{2}\\
c^{2}\end{array}\right)=\frac{1}{\sigma}\left(\begin{array}{c}
\sin^2 \pi\alpha\\
\sin^2 \pi\beta\\
\sin^2 \pi\gamma \end{array}\right),\]
with $\sigma=\sin^2 \pi\alpha+\sin^2 \pi\beta+\sin^2 \pi\gamma$.
With some calculus we can show that\[
ds=\pi J
\left(\begin{array}{c}
d\alpha\\
d\beta\\
d\gamma\end{array}\right),\]
with $J=
\left(\frac{\mbox{diag}(p)}{\sigma}-\frac{sp\transp }{\sigma^{2}}\right)$
and
$
p=\left(\begin{array}{c}
\sin2\pi\alpha\\
\sin2\pi\beta\\
\sin2\pi\gamma\end{array}\right) .
$
The Jacobian determinant that transforms the  angle space to 
squared edge space is  proportional to  $\det\left( \Delta J  \Delta\transp \right)$.
The remaining step, omitted here, is the Jacobian to the hemisphere.
As the experiment begins with triangles uniform on the hemisphere, it is the inverse Jacobian
that we  use in our plot.

 Figure 12 is a Monte-Carlo experiment.   Figure 13, shows
the  theoretical Jacobian. 
\begin{figure}[H]
\begin{center}
\includegraphics[width=250pt]{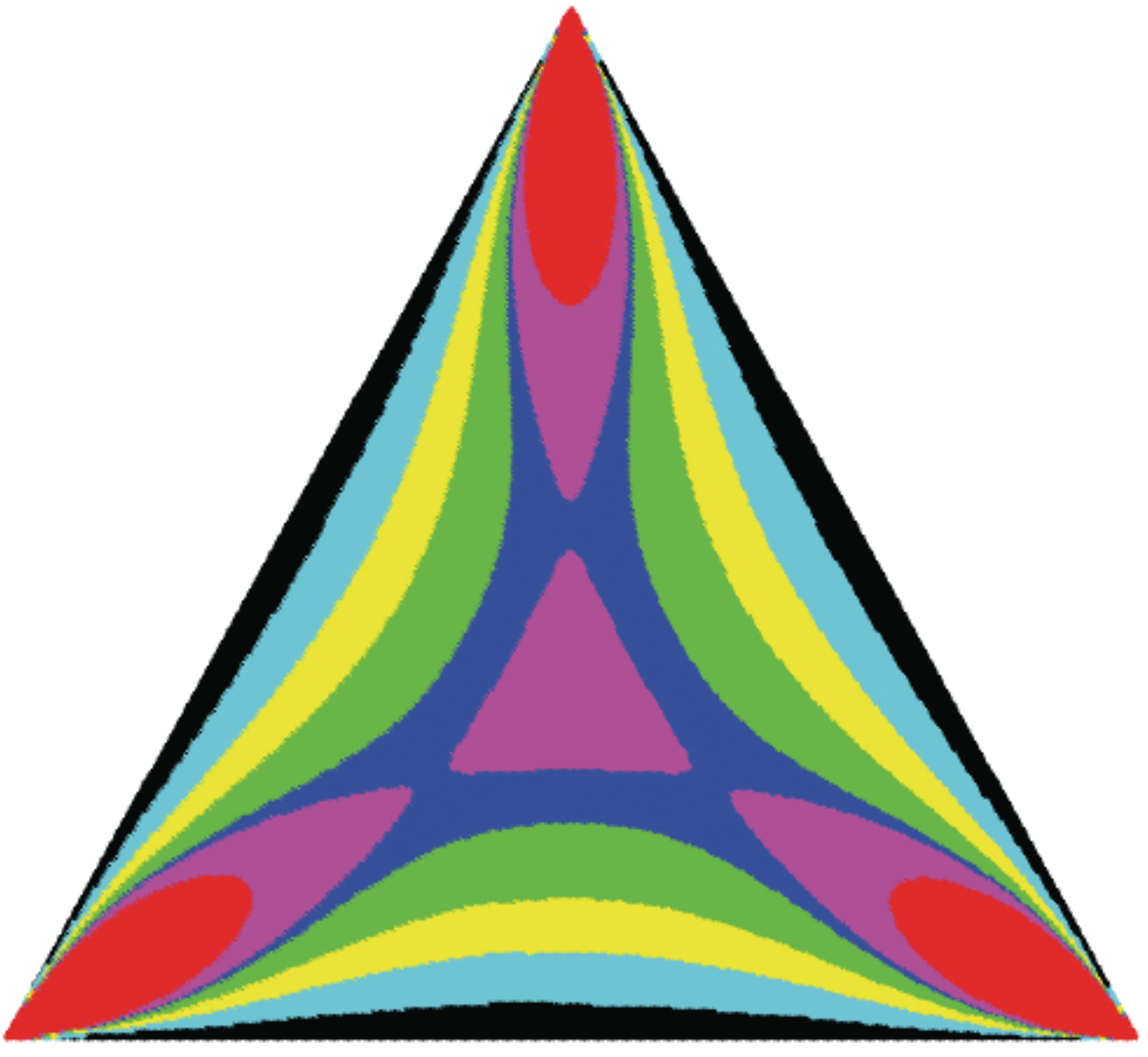}
\end{center}
\caption[]{Angle density  reveals the nonuniformity of angles, when triangles are picked uniformly from shape space. The least likely angles are in black, cyan and yellow with $p<0.1$, $p<0.25,$ and $p<.4$. The most likely angles are in  blue, magenta, and red with $p>0.6$, $p>0.75,$ and $p>0.90$. Green is $.4<p<.6$.}
\label{figure13}
\end{figure}

\section{Further Applications}
A corpus of pictures becomes a set of shapes, and we study a particular feature. A random choice then becomes a random matrix and the entirety of random matrix theory can be considered for shape applications.

We believe that the timing is right to realize the dream of shape
theory. In our favor are 
\begin{itemize}
\item New theorems in random matrix theory that can immediately apply to shape theory
\item Modern computational power making shapes accessible (and  Monte Carlo simulations)
\item The technology of multivariate statistical theory. We can compute hypergeometric functions of matrix argument and zonal polynomials \cite{koev06}. Muirhead \cite{Muirhead82a}  was $30$ years ahead of his time. 
\end{itemize}
There is much room for research in this area. In two dimensions,  the exact condition number distribution for the sphericity test confirms that the shapes come from the standard normal distribution. 

It is very possible that uniformity may not be a great measure for
real problems on shape space. Similar issues arise for random matrices
and random graphs. It is a mistake to think of a random object
as being  {}``any old object.'' Usually random objects
have special properties of their own, much like a special matrix,
or graph, or shape. 

Here are random shapes and convex hulls.
The general technique is multiplying a matrix of standard normals by a $\Delta_n$ from the Helmert matrix.
 It is a bit further afield (but interesting) to ask for the average number of edges of the  convex hull. For literature in this direction see 
\cite[Chapter 8]{schneider08}.
\begin{figure}[H]
\begin{center}
$
\begin{array}{l}
k=3:\\[.73in]
k=5:\\[.73in]
k=20:\\[.73in]
k=100: \\ [3.2in]
\end{array}
\ $
\includegraphics[scale=1.2]{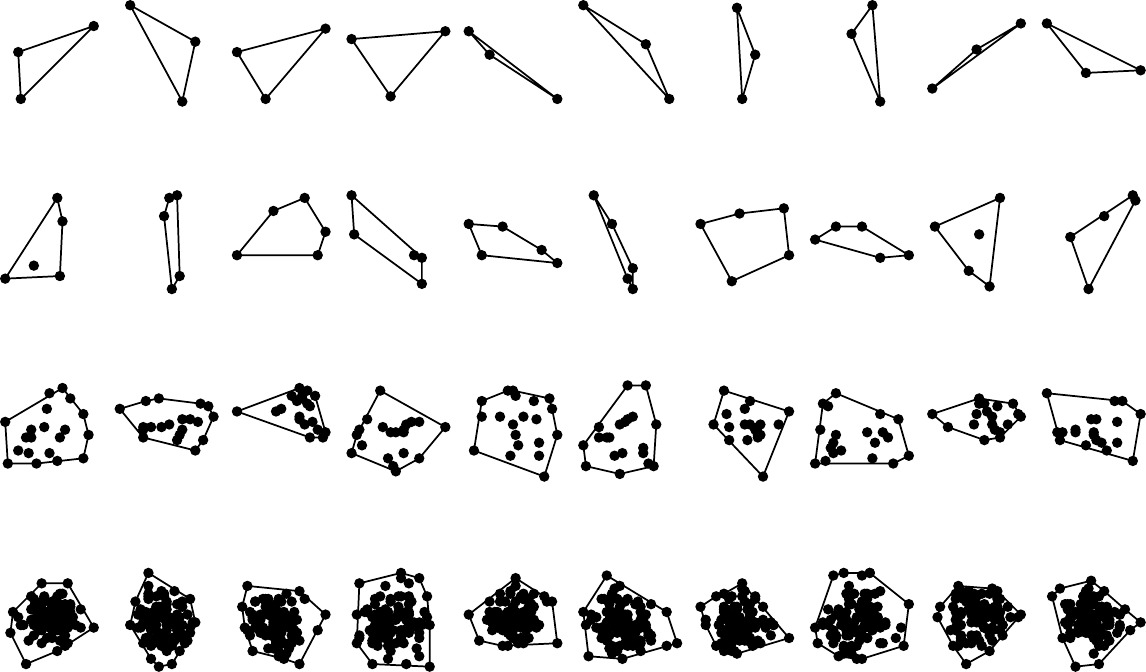}
\end{center}
\vspace{-3in}
\caption[]{Ten random shapes in $2$d taken from the uniform shape distribution with $k=3,5,20,100$ points.}
\end{figure}
\vspace{-0.17in}
\begin{figure}[H]
\begin{center}
\includegraphics[scale=.74]{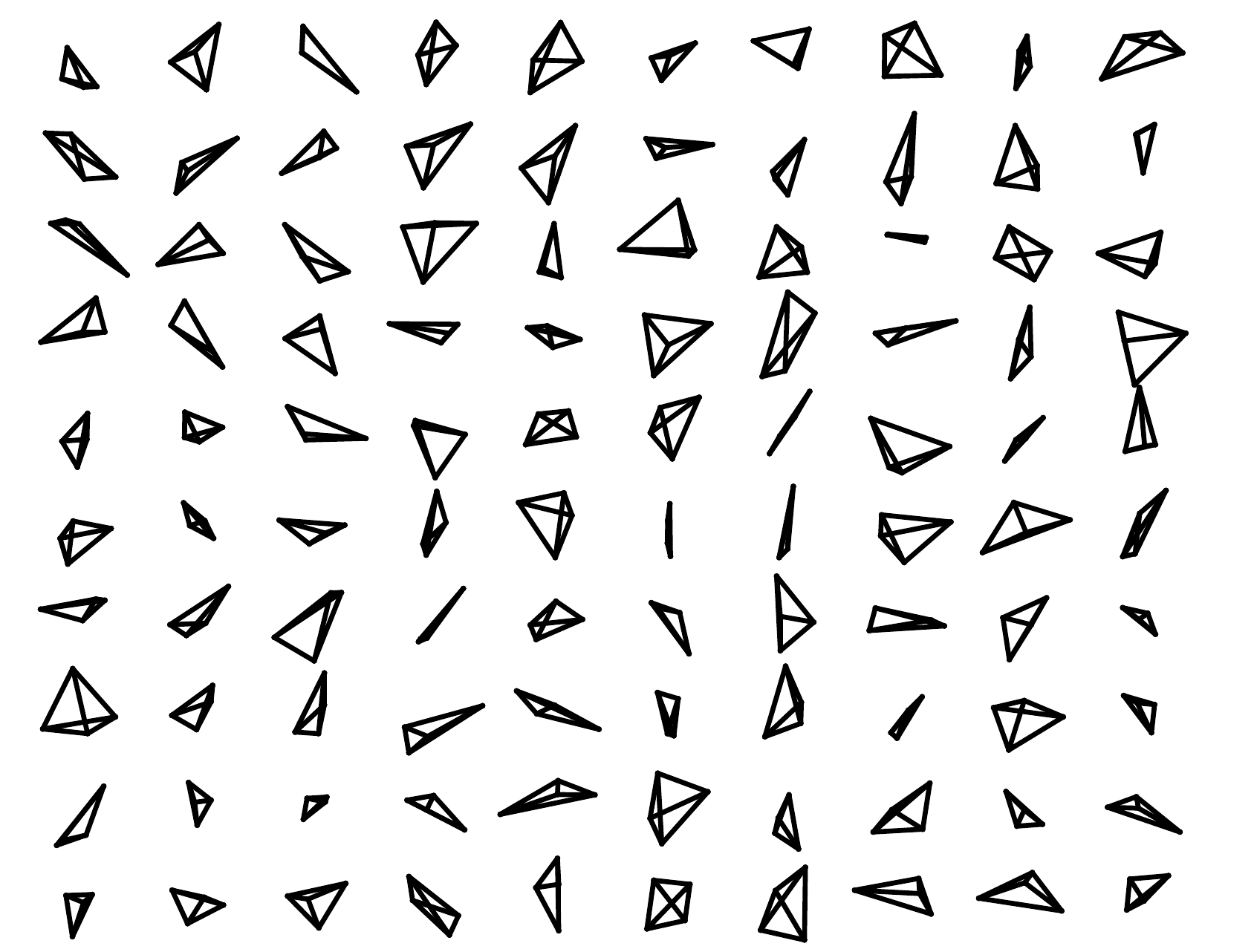}	
\end{center}
\vspace{-0.3in}
\caption[]{ Random tetrahedra in $3$d ($m=3$ and $k=4$)}
\end{figure}
\vspace{-0.132in}
\begin{figure}[H]
\begin{center}
\includegraphics[scale=.95]{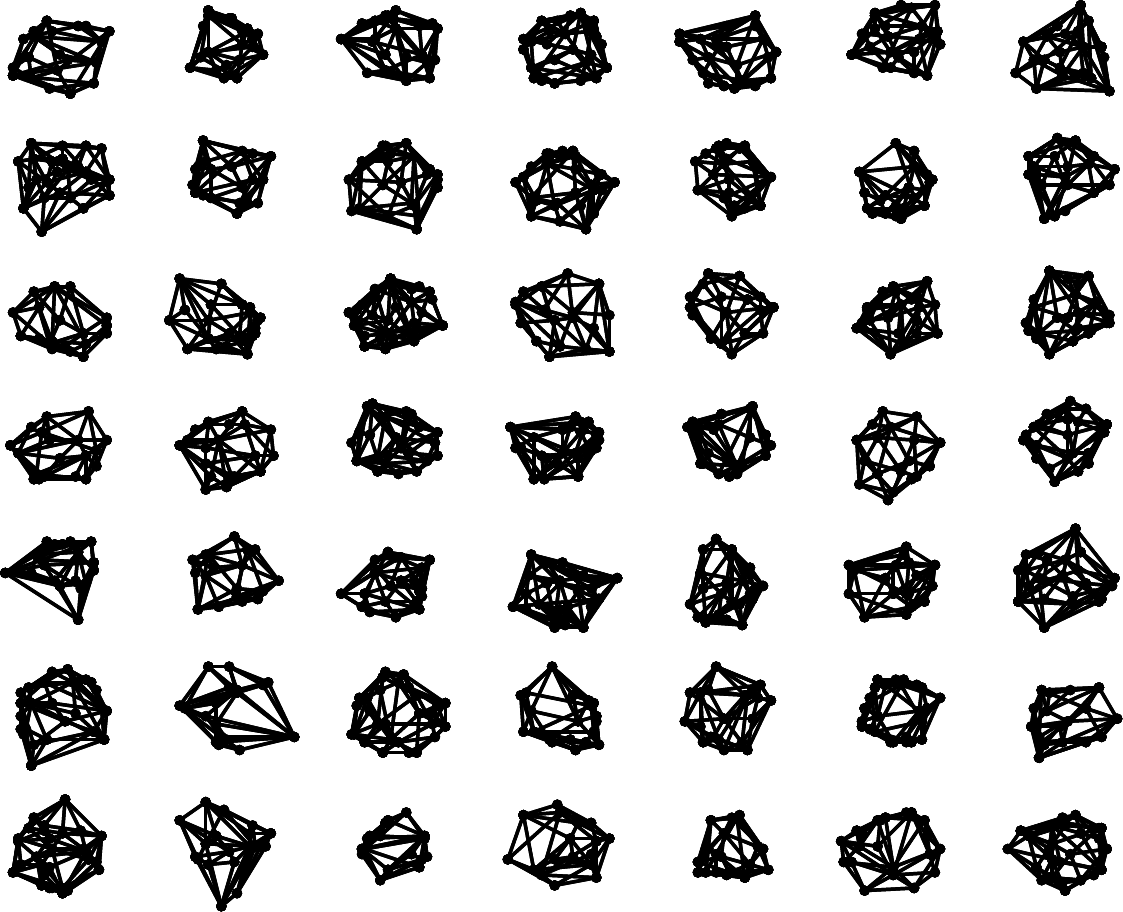}	
\end{center}
\vspace{-0.32in}
\caption[]{Random Gems: Convex hulls of 100 random points in $3$d ($m=3$ and $k=100$)}
\label{fig:gems}
\end{figure}

\subsection{Tests for Uniformity of Shape Space}
An object formed from $k$ points in $\mathbb{R}^{m}$ may be encoded in
an $m\times k$ matrix $X$. The shape is therefore encoded in an $m\times(k-1)$ \emph{preshape matrix} $Z=X\Delta_{k}\transp /\|X\Delta_{k}\transp \|_F,$ where $\Delta_{k}$ is the Helmert matrix in equation (\ref{eq:helmert}). The uniform distribution on shapes may be thought of as the distribution obtained when
\[
X=\mbox{\texttt{\textbf{randn(m,k-1)}}}*\Delta_{k},
\]
so that $X$ is the product of an $m\times(k-1)$ matrix
of iid standard Gaussians and $\Delta_{k}$. In the preshape, $Z$ is normalized by its own Frobenius norm, $\sqrt{\sum Z_{ij}^2},$:
\[
Z=\mbox{\texttt{\textbf{randn(m,k-1)}}}/\|\cdot\|_{F}.
\]
Figure 15 plots random tetrahedra in this way.

In \cite{chikuse04}, Chikuse and Jupp propose a statistical test
on $Z$ for uniformity. From samples $Z_{1},\ldots,Z_{t}$ they
calculate
\[
S=\frac{(k-1)(m(k-1)+2)}{2}\,t\,\mbox{ trace(\ensuremath{\left\{ \ensuremath{\frac{1}{t}\sum_{i=1}\transp Z\transp Z-\frac{1}{k-1}I_{k-1}}\right\} ^{2})\mbox{\ensuremath{\mbox{\ensuremath{}}}}}}
\]
As $t\rightarrow\infty$ they approximate $S$ with $\chi_{(k-1)(k+2)/2}^{2}$. Corrections are proposed for finite $t$ as well.

Other tests are easy to construct. For example, much is known about the smallest singular value of the random matrix $Z$. When $Z$ is square and $t\geq m$, the density of $1/\sigma_{\min}(Z)$ is exactly
\[
\frac{2m\Gamma(\frac{m+1}{2})\Gamma(\frac{m^{2}}{2})}{\sqrt{\pi}\Gamma(\frac{m(m+1)}{2}-1)}t^{1-m^{2}}(t^{2}-m)^{\frac{m(m+1)}{2}-2}{}_{2}F_{1}\left(\frac{m-1}{2},\frac{m}{2}+1;\frac{m^{2}+m}{2}-1;-(t^{2}-m)\right).
\]
Non-square cases can also be handled by a combination of known techniques. 

With the exact density, one can compute the smallest singular value of the samples and perform goodness-of-fit tests (such as Kolmogorov-Smirnov). Broadly speaking one might choose to test for the orthogonal invariance of the singular vectors.  Since $Z*\chi_{m(k-1)}$is normally distributed, tests based on Wishart matrices are natural.

\subsection{Northern Hemisphere Map}
Perhaps  because we had the technology,
we mapped the northern hemisphere (shape theory hemisphere) to angle space (Figure 17.)
Barycentric coordinates correspond to the angles divided by $\pi$.
The resulting picture appears in the figure that follows. The middle
triangle consists of the {}``acute'' points.
\begin{figure}[H]
\begin{center}
\includegraphics[scale=0.7]{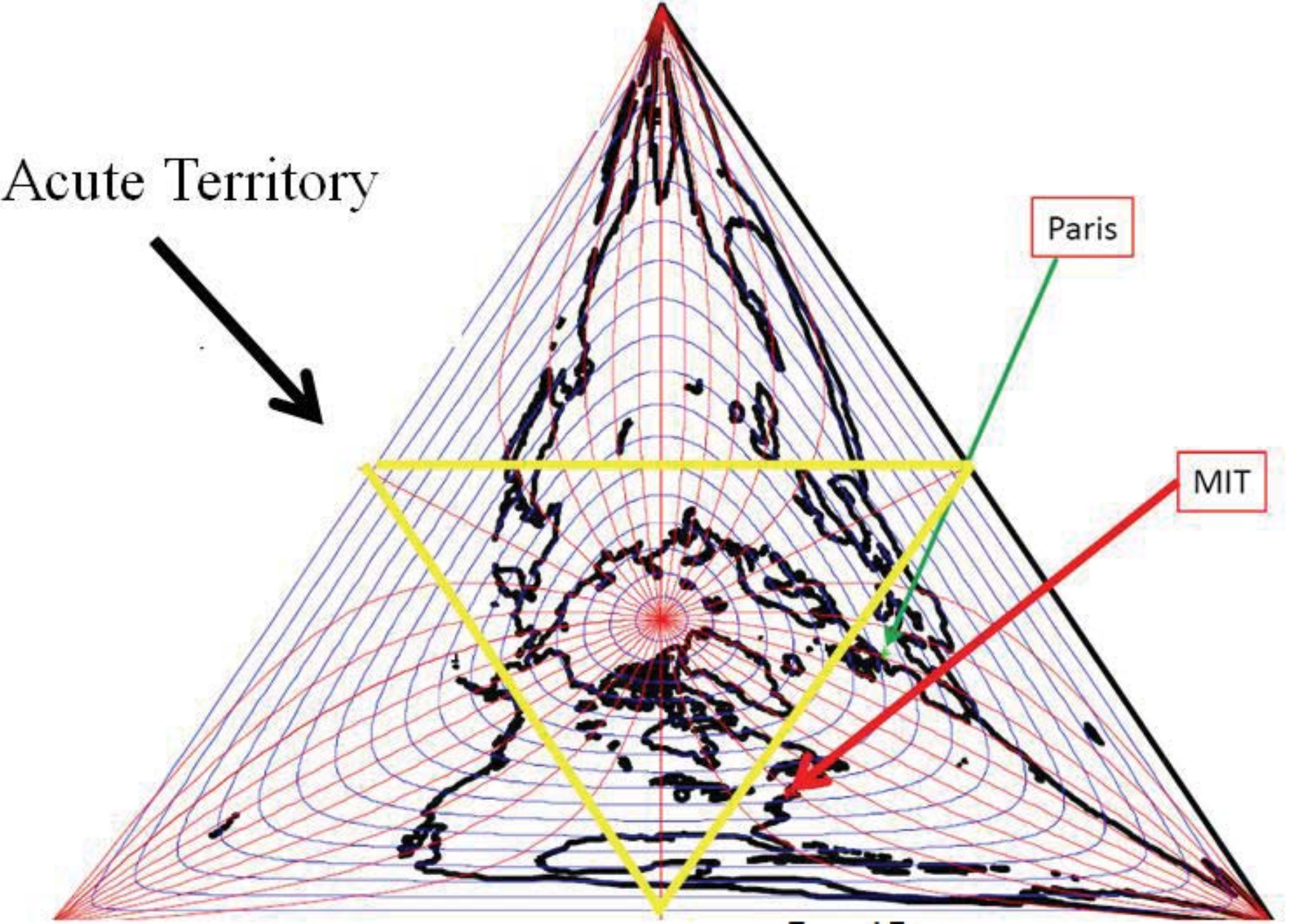}	
\end{center}
\caption[]{Northern hemisphere as triangles in shape space, mapped to {}``angle'' space.}
\end{figure}

We would like to thank Wilfrid Kendall, Mike Todd, and Eric Kostlan for their insights.
The first author acknowledges NSF support under DMS 1035400 and DMS 1016125.
The second author acknowledges NSF support under EFRI 1023152.

\newpage{}
\bibliographystyle{plain}
\bibliography{shape2014.bib}
\end{document}